\documentclass{article}

\usepackage[bbgreekl]{mathbbol}
\usepackage{hyperref}
\usepackage{amsmath}
\usepackage{amssymb}
\usepackage{theorem}
\usepackage{xypic}
\usepackage[all,v2]{xy}
\xyoption{2cell}
\UseAllTwocells
\usepackage{enumitem}
\usepackage{rotating}
\usepackage{mathrsfs}
\usepackage{stmaryrd}
\usepackage{mathdots}
\usepackage[dvipsnames]{xcolor}

\newenvironment{items}
{\begin{enumerate}[topsep=3pt, itemsep=3pt, parsep=0pt, label=(\roman*)]}
{\end{enumerate}}

\newcommand{\tpdf}{\texorpdfstring}

\renewcommand{\tilde}{\widetilde}

\newcommand{\qed}{{\nolinebreak $\,\,\Box$}}

\newcommand{\reg}{{\rm reg}}
\newcommand{\sing}{{\rm sing}}

\newcommand{\A}{{\mathscr A}}

\newcommand{\MM}{{\mathfrak M}}

\newcommand{\Mm}{{\mathfrak m}}

\newcommand{\zz}{{\mathbb Z}}

\newcommand{\pp}{{\mathbb P}}
\newcommand{\cc}{{\mathbb C}}

\newcommand{\ff}{{\mathbb F}}

\renewcommand{\O}{{\mathscr O}}

\newcommand{\iI}{{\mathscr{I}}}

\newcommand{\red}{{\mathop{\rm red}\nolimits}}
\newcommand{\supp}{\mathop{\rm supp}}

\newcommand{\Pic}{\mathop{\rm Pic}\nolimits}

\newcommand{\ol}{\overline}

\newcommand{\dimm}[1]{{^{\color{red} #1}}}

\newcommand{\doublearrowstack}[2]%
                      {{{{\scriptstyle#1}\atop{\textstyle\longrightarrow}}\atop{{\textstyle\longrightarrow}\atop{\scriptstyle#2}}}}
\newcommand{\rightleftarrowstack}[2]%
                      {{{{\scriptstyle#1}\atop{\textstyle\longrightarrow}}\atop{{\textstyle\longleftarrow}\atop{\scriptstyle#2}}}}
\newcommand{\leftrightarrowstack}[2]%
                      {{{{\scriptstyle#1}\atop{\textstyle\longleftarrow}}\atop{{\textstyle\longrightarrow}\atop{\scriptstyle#2}}}}

\newcommand{\sfrac}[2]{{\textstyle\frac{#1}{#2}}}

\newtheorem{thm}{Theorem}[section]
\newtheorem{cor}[thm]{Corollary}
\newtheorem{lem}[thm]{Lemma}
\newtheorem{prop}[thm]{Proposition}

\theorembodyfont{\upshape}
\newtheorem{defn}[thm]{Definition}

\newtheorem{rmk}[thm]{Remark}

\newenvironment{pf}{\begin{trivlist}\item[]{\sc Proof.}}%
            {\nolinebreak $\Box$ \end{trivlist}}

            {\nolinebreak $\Box$ \end{trivlist}}

\newcommand{\EE}{\mathfrak E}

\newcommand{\cok}{\mathop{\rm cok}\nolimits}

\newcommand{\spec}{\mathop{\rm Spec}\nolimits}

\newcommand{\Sym}{\mathop{\rm Sym}\nolimits}

\newcommand{\noprint}[1]{}

\def\Label#1{\label{#1}{\tt [#1]}\phantom{Mh}}
\newcommand{\comment}[1]{$\mbox{}^{\spadesuit}${\marginpar{\footnotesize #1}}}

\def\Label{\label}
\def\comment{\noprint}

\renewcommand{\dimm}[1]{}

\title{On the stability of regular  algebras}
\author{Kai Behrend and Junho Hwang}

\begin{document}
\sloppy

\maketitle
\tableofcontents

\section*{Introduction}
\addcontentsline{toc}{section}{\protect\numberline{}Introduction}%

In \cite{BehrendNoohi}, a notion of stability was introduced for
non-commutative graded algebras (connected and generated in degree~1).
Moreover, quasi-projective moduli 
stacks of Deligne-Mumford type where constructed, for stable graded
algebras of fixed Hilbert series.  Here we are concerned with the
Hilbert series $(1-t)^{-3}$, which is the Hilbert series of the
projective plane. 

Interesting examples of algebras with this Hilbert series are {\em
  quadratic 
  regular algebras }of global dimension~3
(see~\cite{AS},~\cite{ATV}). 
In this article, we determine exactly which of these algebras are
stable, and then describe the moduli stack of stable regular
algebras. 

In \cite{ATV}, it was proved that all quadratic regular algebras of
dimension~3 are
associated to {\em elliptic triples }$(X,\sigma,L)$, consisting of a
scheme $X$, an automorphism $\sigma$ of $X$, and  line bundle $L$,
which embeds $X$ into $\pp^2$ as a divisor of degree~3.  We prove that
the regular algebra $A$ determined by such a triple is stable if and
only if $X$  does not contain any linear component preserved by $\sigma$,
and it does not contain any singularity of $X$ fixed by $\sigma$.  This
leaves only two stable cases: the case where $X$ is non-singular, and
the case where $X$ is a Neron triangle and  the automorphism
$\sigma$ acts transitively on the set of its  components/sides.

See Section~\ref{mainthem} for a more detailed outline of the proof of
this result, which is the main contribution of this paper. 

The moduli stack of these stable algebras turns out to be smooth, and
to decompose into 4 irreducible components, of dimensions 2, 1, 0, and
0, respectively, see Theorem~\ref{stackthm}.  These components  can be
characterized by the order of the automorphism which $\sigma$ induces on
$\Pic^0(X)$.  This order can be 1,2,3, or~4. In the classification of
\cite{AS}, these components correspond to Types A,B,E, and H,
respectively. Unfortunately, neither of the two positive-dimensional
components is proper. 

We prove that the property `stable {\em and }regular' is an open property,
for flat families of graded algebras. Therefore, the moduli stack we
construct is open in the moduli stack of all stable algebras of
Hilbert series $(1-t)^{-3}$, and hence also dense in each component
which it intersects. We do not address the question of whether or not
there are components of the moduli stack of
stable algebras which do not contain any regular algebras.

Let us emphasize that our goal here is not to classify quadratic
regular algebras, or describe their moduli. This has been done
elsewhere, see~\cite{ATV} for the classification, and
\cite{Ueda}, and references therein,  for moduli. Rather, the purpose
of this work is to determine how our notion of stability relates to
these known moduli spaces.

The question of compactifying our moduli stack turns out to be rather
subtle: the geometric invariant theory considerations of
\cite{BehrendNoohi}  provide us with compactifications using
semi-stable algebras, but these compactifications depend on an integer
$q$, because they are constructed in terms of algebras truncated
beyond degree $q$.  Preliminary considerations show 
that the boundaries of these compactifications grow, as $q$ increases.

There is a natural  absolute notion of semi-stability, but it does not
give rise to compact moduli. We conjecture that all quadratic regular
algebras of dimension~3 are  semi-stable. Examples show that the
converse is not true: there are semi-stable algebras which are
degenerate, in the sense of~\cite{ATV}. On the other hand, we do not
know any examples of stable algebras which are not regular. For
details on semi-stability, see~\cite{Junho}.

We will work over an algebraically closed field of characteristic
avoiding~2 and~3, and call it $\cc$. 
All our graded algebras $A$ will be connected, i.e., 
$$A=\bigoplus_{n\geq0}A_n$$
with $A_0=\cc$, and generated in degree~1. 

\section{Stability of regular algebras}

We start by reviewing the notion of stability for graded algebras
from~\cite{BehrendNoohi}, to set up notation.  Then we review the
the theory of quadratic regular algebras of dimension~3, as far as it
is relevant for us. We announce our main result, characterizing the
stable regular algebras, and outline the proof.  

We prove some general facts that we will later need in the proof.

\subsection{Stability}

Let $A$ be a connected graded algebra, finitely generated in
degree~1. We recall some notions from~\cite{BehrendNoohi}. 

\begin{defn}
The algebra $A$ is {\bf $q$-stable}, (for an integer $q>1$), if for
every non-trivial test configuration for $A$, generated in degree 1,
the Futaki function satisfies $F(q)>F(1)$.  It is {\bf stable}, if
there exists an $N>0$, such that it is $q$-stable for all $q>N$.
\end{defn}

The degree 1 part of a test configuration is a filtration 
\begin{equation}\Label{filter}
A_1=W^{(0)}\supset W^{(1)}\supset\ldots\,,
\end{equation}
(the inclusions are not strict), such that $W^{(k)}=0$, for
sufficiently large $k$. 
The {\em test configuration }$B\subset A[t,t^{-1}]$ generated by the
filtration (\ref{filter}) is the $\cc[t]$-algebra defined  by 
\begin{equation}\Label{test}
B_n=\bigoplus_{k}t^{-k}\sum_{k_1+\ldots+k_n=k}W^{(k_1)}\ldots
W^{(k_n)}\,.
\end{equation}
Here, for fixed $n>0$, and $k\in\zz$, the sum is over all $n$-tuples
of non-negative numbers $(k_1,\ldots,k_n)$, such that $k_1+\ldots
k_n=k$. If $k<0$, the set of such $n$-tuples is empty, and the
coefficient of $t^{-k}$ is $A_n=(A_1)^n$, by definition. 

Setting
$$I_n^{(k)}=\sum_{k_1+\ldots+k_n=k}W^{(k_1)}\ldots W^{(k_n)}$$
defines a descending sequence of two-sided graded ideals $A\supset
I^{(1)} \subset I^{(2)}\ldots$ in $A$, satisfying
$I^{(k)}I^{(\ell)}\subset I^{(k+\ell)}$, for all $k,\ell\geq0$.  The
fibre over $t=0$ of $B$ is the doubly graded algebra
$$B/tB=\bigoplus_{k\geq0} I^{(k)}/I^{(k+1)}\,.$$

The {\em weight function }of the test configuration (\ref{test}) is given by 
$$w(n)=\sum_{k>0}
\dim I_n^{(k)}=\sum_{k>0}\dim\Big( \sum_{k_1+\ldots+k_n=k}W^{(k_1)}\ldots
W^{(k_n)}\Big)\,.$$
The weight $w(n)$ is equal to the total weight of the graded vector
space $(B/tB)_n=\bigoplus_{k\geq0}(B/tB)^{(k)}_n$:
$$w(n)=\sum_{k\geq0}k\,\dim (B/tB)^{(k)}_n\,.$$

The {\em Futaki function }of the test configuration $(\ref{test})$ is 
$$F(n)=\frac{w(n)}{n\,\dim A_n}\,.$$
In particular,
$$F(1)=\frac{1}{\dim A_1}\sum_{k>0}\dim W^{(k)}\,.$$

To test for stability, we can restrict to test configurations
generated by filtrations such that $W^{(1)}\not=A_1$.   We will
always  restrict attention to such test configurations. 

\subsection{Quadratic regular algebras of dimension 3}

Recall (\cite{ATV}, Definition~4.5) that an {\bf elliptic triple }is a
triple $(X,\sigma,L)$, consisting of a
scheme $X$, together with a very ample line 
bundle $L$, such that $V=\Gamma(X,L)$ is of dimension~3, and the
closed immersion $X\to
\pp(V)$ makes $X$ into a divisor of degree 3 inside
$\pp(V)$. Moreover, $\sigma:X\to X$ is an automorphism of $X$. 

We will assume our triples to be {\bf regular}, which means that 
$$\sigma^\ast\sigma^\ast L\otimes L\cong \sigma^\ast L\otimes
\sigma^\ast L\,,$$
but not {\bf linear}, which is the stronger condition
$$\sigma^\ast L\cong L\,. $$

The quadratic algebra $A$ associated to the regular elliptic triple
$(X,\sigma,L)$ is the quotient of the tensor algebra of $V$, by the
two-sided ideal generated by the kernel $R$ of
\begin{align}\Label{definerel}
V\otimes V&\longrightarrow \Gamma(X,L\otimes \sigma^\ast L)\\
x\otimes y&\longmapsto x\otimes \sigma^\ast y\,.\nonumber
\end{align}
As (\ref{definerel}) is surjective, and
$\dim\Gamma(X,L\otimes\sigma^\ast L)=6$, the kernel $R$ is of
dimension~3. Therefore, $A$ has 3 generators in degree~1, and 3
relations in degree~2. 

By \cite{ATV} Theorem 6.8.(ii), the algebra $A$ is a {\em regular
  algebra of global dimension~3}, in particular, its Hilbert series
is given by $(1-t)^{-3}$, equivalently, $\dim
A_n=\frac{1}{2}(n+1)(n+2)$, for all $n\geq0$.

The {\em twisted coordinate ring }$B$, associated to the triple
$(X,\sigma,L)$ is defined by 
$$B_n=\Gamma(X,L_n)
\,,$$
where 
$$L_n=L\otimes \sigma^\ast L\otimes
\ldots\otimes(\sigma^{n-1})^\ast L\,,$$ 
which multiplication given by
$$(x_0\otimes\ldots \otimes x_{n-1})\cdot
(y_0\otimes\ldots\otimes y_{m-1})
=x_0\otimes\ldots\otimes x_{n-1}\otimes (\sigma^n)^\ast y_0\otimes
\ldots\otimes(\sigma^n)^\ast y_{m-1}\,.$$
We have $\dim B_n=3n$, for $n>0$. 

There is a canonical  morphism  of graded $\cc$-algebras
\begin{equation}\Label{theemi}
A\longrightarrow B\,,
\end{equation}
induced by the identification $V=A_1=B_1$. By \cite{ATV}
Theorem~6.8, the morphism (\ref{theemi}) is an epimorphism, whose
kernel is equal to both $c_3A$ and $Ac_3$, for an element $c_3\in A_3$,
which is both a left and a right non-zero divisor (of course, $c_3$ is
unique up to multiplication by a non-zero scalar).

\subsection{The main theorem}\Label{mainthem}

\begin{thm}[Stability of regular  algebras]\Label{mnaefo}
Let $(X,\sigma,L)$ be a regular elliptic triple, with associated
regular algebra $A$.  The following are equivalent:
\begin{items}
\item $A$ is stable, 
\item $A$ is $q$-stable, for every $q\geq3$,
\item $A$ is 3-stable,
\item\Label{thmiv} $X$ contains no linear component preserved by
  $\sigma$, and no  
singularity preserved by $\sigma$. 
\end{items}
\end{thm}
(In \ref{thmiv}, it is not required that the linear
component be fixed pointwise by $\sigma$.)

An elliptic triple without linear components preserved by $\sigma$,
and without singularities preserved by $\sigma$, is necessarily
either smooth,
or the union of a line and a smooth conic intersecting at two points,
or a union of three lines intersecting transversally at three nodes.
The second of these three cases is not regular, but {\em exceptional},
see~\cite{ATV},~4.9, and is therefore excluded in the theorem.

Thus, if $X$ is smooth, or a triangle on which $\sigma$ acts by
cyclic permutation of the edges, $A$ is stable, in all other cases $A$
is at best  semi-stable.\comment{forward reference to definition}

\subsubsection{Method of proof}

To prove the theorem, we prove
the following three propositions.

\begin{prop}\Label{sspropU}
Let $U\subset V$ be a one-dimensional subspace, such that the line
$Y=Z(U)\subset \pp(V)$ is contained in $X\subset \pp(V)$, and such that
$\sigma|_{Y}:Y\to 
X$ factors through $Y\subset X$. Then the test configuration for
$A$ generated by the flag $V\supset U\supset 0$ of $A_1$ has constant
Futaki function $\tilde F$, i.e.,  $\tilde F(n)=\tilde F(1)$, for all $n\geq2$. 
\end{prop}

\begin{prop}\Label{sspropW}
Let $W\subset V$ be a 2-dimensional subspace, such that the point
$P=Z(W)\in \pp(V)$ is contained in $X$, is a singular point of
$X$, and satisfies $\sigma(P)=P$.  Then the test configuration for
$A$ generated by the flag $V\supset W\supset 0$ of $A_1$ has constant
Futaki function. 
\end{prop}

\begin{prop}[Stability estimates]\Label{estimates}
Consider the filtration 
\begin{equation}\Label{genfil}
V\supset\underbrace{ W\supset\ldots\supset W}_{\text{$\ell$ times}}\supset
\underbrace{U\supset\ldots\supset U}_{\text{$m$ times}}\supset 0
\end{equation}
of $A_1$, where at least one of the two integers $\ell\geq0$, $m\geq0$
is positive, and $\dim U=1$, $\dim W=2$. Let $P=Z(W)$, and $Y=Z(U)$. 

Suppose that if $P\in X$, and $\sigma(P)=P$, then $P$ is a
non-singular point of $X$, and assume that if $Y\subset X$, then
$X$ is a Neron triangle, whose sides are permuted cyclically by
$\sigma$. Then for the Futaki function $F(n)$ of the test
configuration generated by the flag (\ref{genfil}), we have 
$$
\begin{cases}F(n)\geq F(1)\,, &\text{for $n=1,2$}\\
F(n)>F(1)\,,&\text{for $n\geq 3$}\,.\end{cases}$$
\end{prop}

Propositions \ref{sspropU} and \ref{sspropW} will be proved in the
next section on twisted coordinate rings.

\subsubsection{The stability estimates}

Let us set up notation used in the proof of
Proposition~\ref{estimates}. Let $(X,\sigma,L)$ be a regular elliptic
triple, with associated regular algebra $A$ and twisted coordinate
ring $B=A/c_3 A$. The flag (\ref{genfil}) generates test
configurations $(I_n^{(k)})$ for $B$ and $(J_n^{(k)})$ for $A$.  Let
$w(n)$ be the weight function of $(I_n^{(k)})$ and $\tilde w(n)$ the weight
function of $(J_n^{(k)})$. Let $a(n)=\dim A_n$, and $b(n)=\dim
B_n$, so 
$$a(n)=\frac{1}{2}(n+1)(n+2)\,,\qquad\text{and}\qquad b(n)=3n\,.$$
For future reference, let us also introduce
$$\tilde a(n)=a(n-3)=\frac{1}{2}(n-1)(n-2)\,.$$
Let $F(n)$ be the Futaki function of $(I_n^{(k)})$ and $\tilde F(n)$
the Futaki function of $(J_n^{(k)})$. 

We have 
$$\tilde F(n)=\frac{\tilde w(n)}{\frac{1}{2}n(n+1)(n+2)}\,,$$
and
$$\tilde F(1)=F(1)=\frac{2\ell+m}{3}\,,$$
becase $\tilde w(1)=w(1)=2\ell+m$. 

For $n=2$, we have $\tilde w(2)= w(2)$, and our claim that 
$\tilde F(2)\geq \tilde F(1)$, amounts to 
$$w(2)\geq 8\ell+4m\,.$$

For $n\geq 3$ our claim that $\tilde F(n)>\tilde F(1)$ amounts to 
\begin{equation}\Label{mainclaim}
\tilde w(n)>\displaystyle\frac{2\ell+m}{6}\,n(n+1)(n+2)\,.
\end{equation}
To prove (\ref{mainclaim}), we first prove that lower estimates for
$w(n)$, imply lower estimates for $\tilde w(n)$ via a bootstrapping
method using the short exact sequence
$$\xymatrix{
0\rto& A_{n-3}\rto^-{c_3\,\cdot}& A_n\rto & B_n\rto & 0}\,.$$ 
The larger the integer $p$ such that $c_3\in J^{(p)}_3$, the better the
estimates are that the bootstrapping method gives.  For different
cases, different values for $p$ apply. In many cases, 
it will be sufficient that $c_3\in J^{(2\ell)}_3$, but in the case
where $Y$ is a linear component of $X$, we will need $c_3\in
J_3^{(3\ell)}$. We will prove these facts in Theorem~\ref{anap} and
Proposition~\ref{amplify}.  

Finally, the estimates for $w(n)$ are postponed to the next section
under the heading `Estimates'.

Our estimates for $w(n)$ and $\tilde w(n)$ will always be of the form
\begin{equation}\Label{oftheform}
w(n)\geq \ell\,w_\ell(n)+m\,w_m(n)\qquad\text{and}\qquad
\tilde w(n)\geq\ell\,\tilde w_\ell(n)+m\,\tilde w_m(n)\,,
\end{equation}
with functions $w_\ell(n), w_m(n),\tilde w_\ell(n), \tilde w_m(n)$,
which do not depend on $\ell$ or $m$. 

The stability condition is
$$3\tilde w(n)-(2\ell+m)\,n\,a(n)>0\,,$$
which we can deduce from 
$$\ell\Big(3\tilde w_\ell(n)-2n\,a(n)\Big)
+m\Big(3\tilde w_m(n)-n\,a(n)\Big)>0\,.$$
Let us therefore introduce notation 
$$G_\ell(n)=3\tilde w_\ell(n)-2n \,a(n)=3\tilde
w_\ell(n)-n(n+1)(n+2)\,,$$
and
$$G_m(n)=3\tilde w_m(n)-n \,a(n)=3\tilde
w_m(n)-\sfrac{1}{2}n(n+1)(n+2)\,.$$
With this notation, stability will be implied by
\begin{equation}\Label{ggest}
\ell\, G_\ell(n)+m\,G_m(n)>0\,.
\end{equation}

\begin{defn}\comment{may get rid of this?}
The one-dimensional subspace $U\subset V$ corresponding  to the  line
$Z=Z(U)\subset \pp(V)$ is called  {\bf special}, if $Z$ is contained
in $X$, i.e., if $Z$ is a component of $X$.
\end{defn}

In the non-special case, the zero locus of $U$ in $X$ is an effective
Cartier divisor $D$ on $X$, such that $L=\O(D)$. 

\begin{defn}
The two-dimensional subspace $W\subset V$ corresponding  to the point
$P=Z(W)\in \pp(V)$ is called  {\bf special}, if $P$ lies
on the smooth part of $X$ and $\sigma(P)=P$.
\end{defn}

One of the most useful facts about the non-special case is that it
implies $VW+WV=B_2$, as we will see later. 

If $W$ is special, then $P$ is an effective Weil divisor on $X$, and
we have $W=\Gamma\big(L(-P)\big)$. 

\subsubsection{Bootstrapping: comparing \texorpdfstring{$A$}{A} and
  \texorpdfstring{$B$}{B} using \texorpdfstring{$c_3$}{c3}}

\begin{prop}\Label{crom}
Suppose that $c_3\in J^{(p)}_3$. Then we have
$$\tilde w(n)\geq \sum_{1\leq n-3i\leq n}\big(w(n-3i)+p\,
a(n-3i-3)\big)\,,$$
where the sum is over all integers $i$ satisfying $1\leq n-3i\leq n$. 
\end{prop}
\begin{pf}
By properties of test configurations, we have
$c_3J_{n-3}^{(k)}\subset J_{n}^{(k+p)}$, for all $k\geq0$ and all
$n\in\zz$. Moreover, we have a surjection $J_n^{(k+p)}\twoheadrightarrow
I_n^{(k+p)}$, whose kernel contains $c_3J_{n-3}^{(k)}$, which proves that 
$$\dim J_n^{(k+p)}\geq \dim I_n^{(k+p)}+\dim c_3 J^{(k)}_{n-3}\,,$$
for all $k\geq0$ and all $n\in\zz$. As $c_3$ is a non-zero divisor, we
have
$\dim c_3J_{n-3}^{(k)}=\dim J_{n-3}^{(k)}$, and so 
$$\dim J_n^{(k+p)}\geq \dim I_n^{(k+p)}+\dim  J^{(k)}_{n-3}\,.$$
Thus we have
\begin{align*}
\tilde w(n)&=\sum_{k=1}^\infty \dim J_n^{(k)}\\
&\geq\sum_{k=1}^\infty \dim I_n^{(k)}+p\, \dim A_{n-3}+\sum_{k=1}^\infty \dim
J_{n-3}^{(k)}\\
&=w(n)+p\, a(n-3)+\tilde w(n-3)\,,
\end{align*}
because for $k\leq0$, we have $J_{n-3}^{(k)}=A_{n-3}$. From this
recursion, we deduce
$$
\tilde w(n)\geq\sum_{i=0}^\infty w(n-3i)+p\sum_{i=1}^\infty
a(n-3i)\,,$$
which is our claim.
\end{pf}

\begin{cor}
Assume that  $c_3\in J_3^{(q\ell+sm)}$. Then 
$$\tilde w(n)\geq \sum_{i=0}^{\lfloor\frac{n-1}{3}\rfloor}\big(w+(q\ell+sm)\tilde
a\big)(n-3i)\,,$$
and hence we may use 
$$\tilde w_\ell(n)=
\sum_{i=0}^{\lfloor\frac{n-1}{3}\rfloor}
(w_\ell+q\tilde a)(n-3i)$$
and
$$\tilde w_m(n)=
\sum_{i=0}^{\lfloor\frac{n-1}{3}\rfloor}
(w_m+s\tilde a)(n-3i)\,.$$
\end{cor}

Both $w_\ell(n)$ and $w_m(n)$ will turn out to be polynomial functions
of degree 2 in $n$.  So we will need the following sums later on.

\begin{lem}\Label{sigmanottt}
\begin{align*}
\zeta_2(n)
&=3\sum_{i=0}^{\lfloor\frac{n-1}{3}\rfloor}(n-3i)^2
=\frac{1}{6}\,n(n+3)(2n+3)+
\begin{cases} 
0& \text{if $n\equiv 0\mod 3$}\\
-\frac{1}{3}& \text{if $n\equiv 1\mod 3$}\\
+\frac{1}{3}& \text{if $n\equiv2 \mod 3$\rlap{\,,}}\end{cases}\\
\zeta_1(n)
&=3\sum_{i=0}^{\lfloor\frac{n-1}{3}\rfloor}(n-3i)
=\frac{1}{2}\,n(n+3)+
\begin{cases} 
0& \text{if $n\equiv 0\mod 3$}\\
+1& \text{if $n\equiv 1\mod 3$}\\
+1& \text{if $n\equiv2 \mod 3$\rlap{\,,}}\end{cases}\\
\zeta_0(n)
&=3\sum_{i=0}^{\lfloor\frac{n-1}{3}\rfloor}1
=n+
\begin{cases} 
0& \text{if $n\equiv 0\mod 3$}\\
+2 & \text{if $n\equiv 1\mod 3$}\\
+1& \text{if $n\equiv2 \mod 3$\rlap{\,.}}\end{cases}\\
\end{align*}
\end{lem}

\begin{rmk}
Note that, no matter the divisibility of $n$ modulo $3$,
\begin{equation}\Label{amazing}
n(n+1)(n+2)=3\,\zeta_2(n)-3\,\zeta_1(n)+2\,\zeta_0(n)\,.
\end{equation}
Because of this, it will turn out to be convenient to express all our
estimates for $G_\ell(n)$ and $G_m(n)$ in terms of $\zeta_2(n)$,
$\zeta_1(n)$ and $\zeta_0(n)$. 
\end{rmk}

\subsection{Twisted homogeneous coordinate rings}

Suppose $(X,\sigma,L)$ is an elliptic triple, with associated twisted
homogeneous  coordinate ring 
$$B=\bigoplus_n \Gamma(X,L_n)\,,\qquad L_n=
L\otimes L^\sigma\otimes\ldots\otimes
L^{\sigma^{n-1}}\,.$$

Let $(Z,\tau)$ be another scheme and automorphism, and let $\phi:Z\to
X$ be a morphism, such that $\sigma\circ\phi=\phi\circ\tau$. Let
$N=\phi^\ast L$. Let 
$$\tilde B=\bigoplus_n\Gamma(Z,N_n)\,,\qquad N_n=N\otimes
N^\tau\otimes\ldots\otimes 
L^{\tau^{n-1}}$$
be the twisted homogeneous coordinate ring of $(Z,\tau,N)$. 

We have a canonical ring morphism 
$$B\to \tilde B\,.$$
In degree $n$ it is given by the pullback map
$$(\phi^\ast)^{\otimes n}:\Gamma(X,L_n)\longrightarrow 
\Gamma(Z,N_n)\,,$$
which exists because 
$$N^{\tau^i}=(\phi^\ast L)^{\tau^i}=\phi^\ast(L^{\sigma^i})\,.$$

\subsubsection{Line in \texorpdfstring{$X$}{X}}

Let $U\subset V$ be 1-dimensional, and $Z=Z(U)\subset \pp(V)$. Assume
that $Z\subset X$, and that $\sigma$ factors through $Z$. Then we get
an associated algebra quotient $B\to \tilde B$. 
The algebra $\tilde B$ is a twist of the polynomial ring in two
variables, or a quantum projective line. 

In degree 2, we have a short
exact sequence
$$\xymatrix{
0\rto & Q\rto & V/U\otimes V/U \rto & \tilde B_2\rto & 0}\,.$$
The dimension of $\tilde B_2$ is 3, by the classification of quantum
projective lines, and $\dim Q=1$.  We consider the induced morphism of
short exact sequences
$$\xymatrix{
0\rto & R\rto\dto & V\otimes V \rto\dto & B_2\rto\dto & 0\\
0\rto & Q\rto & V/U\otimes V/U \rto & \tilde B_2\rto & 0\rlap{\,.}}$$

Let $T(V)$ be the free or tensor algebra on $V$, and $A=T(V)/R$ the
regular algebra associated to our elliptic triple. Let $C_3\subset
A_3$ be the one-dimensional subspace generated by $c_3$.

\begin{lem} The map $R\to Q$ is surjective.\end{lem}
\begin{pf}
If not, it is zero. Then $R$ is in the kernel of $V\otimes V\to
V/U\otimes V/U$, which is $U\otimes V+V\otimes U$. If we choose a
basis $(x_i)$ for $V$, with $x_1\in U$, and a basis $(f_i)$ for $R$,
then the matrix $M$, such that $f=Mx$ is of the form
$$M=\begin{pmatrix}
m_{11} & a_{12}\,x_1 & a_{13}\,x_1 \\
m_{21} & a_{22}\,x_1 & a_{23}\,x_1 \\
m_{31} & a_{32}\,x_1 & a_{33}\,x_1 \end{pmatrix}\,.
$$
Because $X$ is defined by $\det M=0$, we see that the double line
$x_1^2=0$ is contained in $X$. Also, we see that along this line, the
matrix $M$ has rank at most 1.  This contradicts the assumption that
$A$ is regular, by Lemma~4.4 of~\cite{ATV}. 
\end{pf}

\begin{cor}
We have $\tilde B=T(V)/\langle U,R\rangle$,
where $\langle U,R\rangle$ denotes the two-sided ideal in the tensor
algebra $T(V)$, generated by the subspaces $U\subset  V$ and $R\subset
V\otimes V$.  
\end{cor}
\begin{pf}
By construction, $\tilde B$ comes with a morphism $T(V)/\langle U,R\rangle\to
\tilde B$.  We have to construct a morphism in the other direction. 

By the classification of quantum projective lines, we know that
$\tilde B=T(V/U)/\langle Q\rangle$.  
The embedding $V/U\to T(V)/\langle U\rangle $ induces an algebra
morphism $T(V/U)\to T(V)/\langle U\rangle$, hence a 
morphism $T(V/U)\to T(V)/\langle U,R\rangle$, which, by the lemma,
annihilates $Q$. This gives us the morphism $\tilde B\to T/\langle
U,R\rangle$, which is the required inverse. 
\end{pf}

\begin{cor}\Label{c3line}
In $A_3$, we have $C_3\subset UVV+VUV+VVU$. 
\end{cor}
\begin{pf}
We have the succession of quotients
$$\xymatrix{
T(V)\ar@{->>}[r] &A\ar@{->>}[r] &B\ar@{->>}[r] &\tilde B}\,,$$
which corresponds to the sequence of ideals in $T(V)$
$$
0\subset \langle R\rangle \subset \langle R, C_3\rangle\subset 
\langle U,R\rangle\,.$$
In particular, modulo $\langle R\rangle$, we have $C_3\subset 
\langle U\rangle_3=UVV+VUV+VVU$. 
\end{pf}

\paragraph{Proof of Proposition~\ref{sspropU}.}

In $B$, we have that $UV=VU$. It follows that this is also true in
$A$, as it is a claim about $A_2=B_2$.  Then it follows that
$UA=AU$, because $A$ is generated by $V$. From this, it follows that
in fact $C_3\subset UVV$ and 
$C_3\subset VVU$. Let $z$ be a generator of $U$. Then it follows that
$z$ is a left and right regular element in $A$, because this is true for any
generator $c_3$ of $C_3$. 

The test configuration in $A$, generated by the flag $V\supset
U\supset 0$ in $A_1$, is the sequence of two-sided ideals
$J^{(k)}=\langle U\rangle^k$.  By what we just proved, we have
$\langle U\rangle ^k= z^k A\cong A(-k)$ the shift of $A$ by $-k$. 

Therefore, for the weight function $\tilde w(n)$ of the test
configuration $(J^{(k)})$, we have
$$\tilde w(n)=\sum_{k=1}^n a(k-1)\,,$$
where $a(n)=\dim A_n=\frac{1}{2}(n+1)(n+2)$. This implies that 
$$\tilde w(n)=\frac{1}{3}\,n\, a(n)\,,$$
or $$\tilde F(n)=\frac{\tilde w(n)}{n\,a(n)}=\frac{1}{3}=\tilde
F(1)\,.$$
So we see that the Futaki function of $(z^k)$ is constant.  This
proves Proposition~\ref{sspropU}. \qed

\subsubsection{Line with embedded point}

Consider now $\pp^1$ with a 2-dimensional embedded point. Call this
scheme 
$\tilde Z$, and $Z=\tilde Z^{\red}=\pp^1$. 
The scheme $\tilde Z$ is embedded into $\pp(V)$ by the choice of a
flag $0\subset U\subset W\subset V$, with $\dim U=1$ and $\dim W=2$,
as the intersection of two 
quadrics $\tilde Z=Z(UW)$. The homogeneous coordinate ring of $\tilde
Z$ is $\cc[V]/(UW)$, which is, in coordinates, $\cc[x,y,z]/(x^2,xy)$. 

\begin{lem}
Any twisted homogeneous coordinate ring of $\tilde Z$ is
quadratic. More precisely, if $\sigma:\tilde Z\to \tilde Z$ is any
scheme automorphism, then the associated twisted homogeneous
coordinate ring is 
$$\tilde B=\frac{T(V)}{\langle U\rangle\langle W\rangle+
\langle W\rangle\langle U\rangle+\langle Q\rangle}\,,$$
where $Q$ is defined by the exact sequence
$$\xymatrix{
0\rto & Q\rto & {\displaystyle\frac{ V\otimes V}{ U\otimes
    W+W\otimes U}} \rto & \tilde 
B_2\rto & 0}\,.$$ 
\end{lem}
\begin{pf}
\comment{
I checked this. Basically, every automorphism of $\tilde Z$ comes from
an automorphism of $V$, respecting the flag $U\subset W\subset V$.}
We omit the proof, which is not difficult, as this result is not used
in the proof of the main theorem.
\end{pf}

 Assume now 
that $\tilde Z\subset X$, and that $\sigma$ factors through $\tilde Z$. This will
be the case, for example, if $X$ is nodal, contains the line $Z(U)\subset
X$, invariant by $\sigma$, and a node $Z(W)\subset X$ fixed by
$\sigma$. 

Then we get
an associated algebra quotient $B\to \tilde B$. In degree 2, we have a short
exact sequence
$$\xymatrix{
0\rto & Q\rto & {\displaystyle\frac{ V\otimes V}{ U\otimes
    W+W\otimes U}} \rto & \tilde 
B_2\rto & 0}\,.$$ 
The dimension of $\tilde B_2$ is 4, so $\dim Q=2$.  We consider the
induced morphism of 
short exact sequences
$$\xymatrix{
0\rto & R\rto\dto & V\otimes V \rto\dto & B_2\rto\dto & 0\\
0\rto & Q\rto & {\displaystyle\frac{ V\otimes V}{ U\otimes
    W+W\otimes U}}\rto & \tilde B_2\rto & 0\rlap{\,.}}$$
\begin{lem} The map $R\to Q$ is surjective.\end{lem}
\begin{pf}
If not, the intersection of $R$ and $UW+WU$ has dimension 2.  If we choose a
basis $(x_i)$ for $V$, with $x_1\in U$, and $x_2\in W$, and a basis
$(f_i)$ for $R$, such that $f_1,f_2\in UW+WU$, 
then the matrix $M$, such that $f=Mx$ is of the form
$$M=\begin{pmatrix}
m_{11}(x_1,x_2) & a_{12}\,x_1 & 0 \\
m_{21}(x_1,x_2) & a_{22}\,x_1 & 0 \\
m_{31} & m_{32}& m_{33} \end{pmatrix}\,.
$$
At the point $\langle 0,0,1\rangle$, this matrix has rank 1, so this
contradicts the non-degeneracy assumption, see Lemma~4.4 of~\cite{ATV}. 
\end{pf}

\begin{cor}
Let $T(V)$ be the free algebra on $V$.  Then 
$$\tilde B=\frac{ T(V)}{\langle U\rangle\langle W\rangle+
\langle W\rangle\langle U\rangle+\langle R\rangle}\,.$$
\end{cor}
\begin{pf}
By the lemma, modulo $\langle W\rangle\langle U\rangle$, we have $Q=R$.
\end{pf}

\begin{cor}\Label{c3wort}
In $A_3$, we have $$C_3\subset UWV+UVW+VUW+WUV+VWU+WVU\,.$$
\end{cor}
\begin{pf}
We have the succession of quotients
$$\xymatrix{
T(V)\ar@{->>}[r] &A\ar@{->>}[r] &B\ar@{->>}[r] &\tilde B}\,,$$
which corresponds to the sequence of ideals in $T(V)$
$$
0\subset \langle R\rangle \subset \langle R, C_3\rangle\subset 
\langle U\rangle\langle W\rangle+
\langle W\rangle\langle U\rangle+\langle R\rangle\,.$$
In particular, modulo $\langle R\rangle$, we have $C_3\subset 
\big(\langle U\rangle\langle W\rangle+
\langle W\rangle\langle U\rangle\big)_3$. 
\end{pf}

\subsubsection{Singular point}

Consider now a point $P=Z(W)\in X$, fixed by $\sigma$, and assume that
$P$ is a singularity of $X$. As $\sigma:X\to X$ is a scheme
automorphism, $\sigma$ induces an automorphism of the first order
neighbourhood $Z$ of $P$ in $X$, which is isomorphic to $\spec
\cc[x,y]/(x,y)^2$, as $P$ is a singular point of $X$. The scheme
$Z=\spec\cc[x,y]/(x,y)^2$ is embedded into $\pp(V)$ as $Z=Z(WW)$. 

Let $\tilde B$ the the twisted coordinate ring of $(Z,\sigma|_Z)$. We
have quotients $A\to B\to \tilde B$. Choose a vector $z\in V$,
$z\not\in W$. Then $\tilde z$, the image of $z$ in $\tilde B$, is a
left and right regular element.  Moreover, both left and right
multiplication by 
$\tilde z$ induce  isomorphisms  $\tilde B_n\to \tilde B_{n+1}$,
for all $n\geq2$.  

Let $Q$ be the kernel defined by the exact sequence
$$\xymatrix{
0\rto & Q\rto & \displaystyle \frac{V\otimes V}{W\otimes W}\rto &
\tilde B_2\rto & 0}\,.$$
Then $\tilde B=T(V)/(\langle WW\rangle+\langle Q\rangle)$, as can be
seen by studying the structure of $\tilde B$. 

As above, regularity implies that  $R$ maps onto $Q$, and hence that 
$$\tilde B=\frac{T(V)}{\langle W\rangle\langle W\rangle +\langle
  R\rangle}\,,$$
and that 
\begin{equation}\Label{crwww}
c_3\in WWV+WVW+VWW\,.
\end{equation}

In particular, we have that left or right multiplication by $z\in A$
induces an isomorphism  
$$\big(A/\langle W\rangle^2\big)_n\stackrel{z}{\longrightarrow} 
\big(A/\langle W\rangle^2)_{n+1}\,,$$
for all $n\geq2$. 

Consider now the (doubly) graded ring 
$$R=\bigoplus_{k=0}^\infty \langle W\rangle^n/\langle
W\rangle^{n+1}\,.$$
In fact, $R$ is the central fibre of the test configuration generated
by $W$, and $A$ is S-equivalent to $R$. 

By the facts we proved, the subring of $k$-degree zero, $R^0=A/\langle
W\rangle$ is a quantum 
projective line, and $R$ is generated over $R^0$ by one element which
quantum commutes with $R^0$.  From this, the graded dimension of each
$J^{(k)}\langle W\rangle^k$ can be computed, giving a proof of
Proposition~\ref{sspropW}. 
\comment{This proof might be fleshed out a bit more.}

\subsection{More analysis of  \texorpdfstring{$c_3$}{c3}}\Label{sec16}

Let $(X,\sigma,L)$ be a regular elliptic triple, such that $L\not\cong
L^\sigma$.

Recall the notion of {\bf tame }line bundle, defined in~\cite{ATV}. 
A line bundle on $X$ is tame if either  $H^0$ vanishes, or $H^1$ vanishes,
or if it is trivial.  For tame line bundles, $h^0$ and $h^1$ can
be calculated as $\deg$ or $-\deg$, respectively. Bundles generated by
sections, duals of tame bundles, and bundles with non-negative degree
on each component are all tame. For example, every non-trivial tame line
bundle of degree 0 has $H^0=H^1=0$. 

As in \cite{ATV}, for a sheaf $F$, generated by global sections, we
denote by $F''$ the kernel of the epimorphism 
$\Gamma(X,F)\otimes_\cc\O_X\to F$. The reason for introducing $F''$ is the
following Proposition~7.17 of~\cite{ATV}.

\begin{prop}\Label{717}
If $M$ is a coherent $\O_X$-module generated by global sections, and
$N$ is locally free with $H^1(X,N)=0$, then the kernel and cokernel of
the multiplication map
$$\Gamma(X,M)\otimes_k\Gamma(X,N)\longrightarrow
\Gamma(X,M\otimes_{\O_X}N)$$
are given by $H^0(X,M''\otimes_{\O_X} N)$ and $H^1(X,M''\otimes_{\O_X}N)$,
respectively. 
\end{prop}

If we apply this to $M=L$ and $N=L^\sigma$, we get the exact sequence
$$\xymatrix{
0\rto & \Gamma(L''\otimes L^\sigma)\rto & \Gamma(L)\otimes
\Gamma(L^\sigma)\rto & \Gamma(L\otimes L^\sigma)\rto & 0}\,.$$
If we identify $A_2$ with $\Gamma(L\otimes L^\sigma)$, this identifies
$\Gamma(L''\otimes L^\sigma)$ with $R$, the kernel of multiplication
$V\otimes V\to A_2$.  By the regularity of $A$, we can write $A_3$ as
the cokernel of $R\otimes A_1\to V\otimes A_2$. This gives the
exactness of the lower row in the following diagram:
\begin{equation}\Label{fir}
\vcenter{\xymatrix{
\Gamma(L''\otimes
L^\sigma)\dimm{3}\otimes\Gamma(L^{\sigma^2})\dimm{3}\rto^\beta\dto^= &
\Gamma(L''\otimes 
L^\sigma\otimes L^{\sigma^2})\dimm{9}\rto\ar@{^{(}->}[d] & C_3\dimm{1}\ar@{^{(}->}[d]\rto &0\\
\Gamma(L''\otimes
L^\sigma)\dimm{3}\otimes\Gamma(L^{\sigma^2})\dimm{3}\rto & \Gamma(L)\dimm{3}\otimes\Gamma(L^\sigma\otimes L^{\sigma^2})\dimm{6}\rto &
A_3\dimm{10}\rto & 0\rlap{\,.}}}\end{equation}
The exactness of the upper row comes from Lemma~7.29 of~\cite{ATV}.

\subsubsection{The case where \texorpdfstring{$P$}{P} 
lies on \texorpdfstring{$X$}{X}}

Consider a 2-dimensional subspace $W\subset\Gamma(L)$, which does {\em
  not }generate $L$.  Let $P\in X$ be the point where $W$ fails to
generate $L$. (Mapping $P$ into $\pp^2$, via the embedding
$X\to\pp^2$, we obtain the point in $\pp^2$, dual to the hyperplane
$W$ in $V$.) We have a surjection of sheaves on $X$
\begin{equation}\Label{ipim}
\xymatrix{W\otimes_\cc\O_X\ar@{->>}[r]& L\otimes_{\O_X} \Mm_P=\Mm_PL}\,.
\end{equation}
If $P$ is a smooth point of $X$, then $P$ is a Cartier divisor, and hence
$\Mm_PL=L(-P)$, but we do not want to make this
assumption. 

We note also that $\Gamma(L\otimes_{\O_X}\Mm_P)=W$, and that $\Mm_PL$
is generated by global sections. We consider $(\Mm_PL)''$:
$$\xymatrix{
0\rto & (\Mm_PL)''\rto & W\otimes \O_X\rto & \Mm_P\,L\rto & 0}$$
If $P$ is a smooth point of $X$, then $(\Mm_PL)''$ is locally free of
rank 1, but not otherwise.

\begin{lem}
Any non-zero element of $\Lambda^2W$ induces an isomorphism 
$$\xymatrix{(\Mm_PL)''\rto^-\sim& (\Mm_PL)^\vee}\,.$$
\end{lem}
\begin{pf}
Let us choose a non-zero element in $\Lambda^2 W$.  As $\Lambda^2 W$
is one-dimensional, this element is a basis element, and also induced
a basis element in the dual space $\Lambda^2 W^\vee$.  Let us denote
this element by $\kappa$, it is a non-degenerate alternating pairing
on $W$, and induces a non-degenerate alternating pairing on the
trivial bundle $W\otimes \O_X$. 
 
We claim that $(\Mm_PL)''\subset W\otimes \O_X$ is isotropic with respect to
$\kappa$. As $\kappa$ takes values in $\O_X$, and $\Mm_P$ contains 
a non-zero divisor of $\O_X$, we
can check this claim after removing $P$ from $X$, where it becomes
trivial, as all sheaves involved become locally free.  Thus, $\kappa$
induces a pairing  
$$\xymatrix{\kappa:(\Mm_PL)''\otimes_{\O_X}\Mm_PL\rto & \O_X}\,.$$
This pairing is non-degenerate, i.e., induces injections $(\Mm_PL)''\to
(\Mm_PL)^\vee$ and $\Mm_PL\to {(\Mm_PL)''}^\vee$.  To prove these facts, it is
again enough to restrict to $X-\{ P\}$, because both $(\Mm_PL)''$ and
$\Mm_PL$ are  submodules of locally free finite rank
$\O_X$-modules. 

By the snake lemma, the 
cokernel of $(\Mm_PL)''\to (\Mm_P L)^\vee$ is equal to the kernel of 
$\Mm_PL\to {(\Mm_PL)''}^\vee$, and hence vanishes. 
\end{pf}

\begin{lem}\Label{notprev}
We have 
$$H^1\big((\Mm_PL)''\otimes L^\sigma\big)=0\,,\qquad\text{and}\qquad
H^1\big((\Mm_PL)''\otimes 
L^\sigma\otimes L^{\sigma^2}\big)=0\,.$$
\end{lem}
\begin{pf}
By the previous lemma, there exists a monomorphism $L^{-1}\to
(\Mm_PL)^\vee\cong (\Mm_PL)''$, whose cokernel is supported over $\{P\}$. Hence
we get induced epimorphisms $H^1(L^{-1}\otimes L^\sigma)\to
H^1\big((\Mm_PL)''\otimes L^\sigma\big)$ and $H^1(L^{-1}\otimes
L^\sigma\otimes L^{\sigma^2})\to H^1\big((\Mm_PL)''\otimes
L^\sigma\otimes L^{\sigma^2}\big)$. Now it suffices to remark that by
Lemma~7.18 of \cite{ATV}, 
both line bundles $L^{-1}\otimes L^\sigma$ and $L^{-1}\otimes
L^\sigma\otimes L^{\sigma^2}$ 
are tame, and so their respective $H^1$ vanishes; in the first case
because $L^{-1}\otimes L^{\sigma}$ is a non-trivial line bundle of degree~0,
and in the second case because $L^{-1}\otimes L^\sigma\otimes L^{\sigma^2}$ is of 
degree~3. 
\end{pf}

\begin{lem}\Label{gammas}
The sheaf $\Mm_P L^\sigma$ is also generated by its global
sections. We have $(\Mm_PL^\sigma)''\cong(\Mm_PL^\sigma)^\vee$, as
well as $H^1\big((\Mm_PL^\sigma)''\otimes L^{\sigma^2}\big)=0$.
\end{lem}
\begin{pf}
The same proofs apply.
\end{pf}

Choosing a basis for $V/W$, we get the exact square
\begin{equation}\Label{square}
\vcenter{\xymatrix{
(\Mm_PL)''\rto\dto & W\otimes \O_X\rto\dto & \Mm_PL\dto\\
L''\rto\dto & \Gamma(L)\otimes \O_X\rto\dto & L\dto\\
\Mm_P\rto & \O_X\rto & \O_P}}\end{equation}
For $N=L^\sigma$ and $N=L^\sigma\otimes L^{\sigma^2}$, we consider the
induced short exact 
sequence
$$\xymatrix{0\rto & (\Mm_PL)''\otimes N\rto & L''\otimes N\rto &
  \Mm_PN\rto & 
  0}\,.$$
In both cases, we get an induced short exact sequence of vector spaces
$$\xymatrix{
0\rto & \Gamma\big((\Mm_PL)''\otimes N\big)\rto & \Gamma(L''\otimes N)\rto &
\Gamma( \Mm_PN)\rto & 0}\,,$$
because $H^1\big((\Mm_PL)''\otimes N\big)=0$, by Lemma~\ref{notprev}.

We obtain the morphism of short exact sequences
\begin{multline*}
\xymatrix{
0\rto & 
\Gamma\big((\Mm_PL)''\otimes L^\sigma\big)\dimm{1}\otimes
\Gamma(L^{\sigma^2})\dimm{3}\rto\dto^\alpha &  
\Gamma(L''\otimes L^\sigma)\dimm{3}\otimes
\Gamma(L^{\sigma^2})\dimm{3}\dto^\beta\rto & \\
0\rto & 
\Gamma\big((\Mm_PL)''\otimes L^\sigma\otimes L^{\sigma^2})\dimm{4}\rto & 
\Gamma(L''\otimes L^\sigma\otimes L^{\sigma^2})\dimm{9}\rto &}
\\
\xymatrix{
\rto &  
\Gamma(\Mm_PL^\sigma)\dimm{2}\otimes
\Gamma(L^{\sigma^2})\dimm{3}\rto\dto^\gamma & 0\\ 
\rto &
\Gamma(\Mm_PL^\sigma\otimes L^{\sigma^2}\big)\dimm{5}\rto & 0\rlap{\,.}}
\end{multline*}
By Lemma~\ref{gammas}, $\gamma$ is
surjective. It follows that $\cok\alpha\to \cok\beta$ is an
epimorphism. Hence we get a morphism of exact
sequences
\begin{equation}\Label{sec}
\vcenter{\xymatrix{
\Gamma\big((\Mm_PL)''\otimes
L^\sigma\big)\dimm{1}\otimes\Gamma(L^{\sigma^2})\dimm{3}\rto^\alpha\ar@{^{(}->}[d] &
\Gamma\big((\Mm_PL)''\otimes 
L^\sigma\otimes L^{\sigma^2}\big)\dimm{4}\rto\ar@{^{(}->}[d] &
\cok\alpha\dimm{1}\ar@{->>}[d]\rto &0\\ 
\Gamma(L''\otimes 
L^\sigma)\dimm{3}\otimes\Gamma(L^{\sigma^2})\dimm{3}\rto^\beta &
\Gamma(L''\otimes L^\sigma\otimes L^{\sigma^2})\dimm{9}\rto & 
C_3\dimm{1}\rto & 0
}}\end{equation}

Composing (\ref{sec}) with (\ref{fir}), we obtain the  morphism of
exact sequences
\begin{equation}\Label{monom}
\vcenter{\xymatrix{
\Gamma((\Mm_PL)''\otimes
L^\sigma)\dimm{1}\otimes\Gamma(L^{\sigma^2})\dimm{3}\rto\ar@{^{(}->}[d] &
\Gamma((\Mm_PL)''\otimes 
L^\sigma\otimes L^{\sigma^2})\dimm{4}\rto\ar@{^{(}->}[d] &
\cok\alpha\dimm{1}\ar@{->}[d]\rto &0\\ 
\Gamma(L''\otimes
L^\sigma)\dimm{3}\otimes\Gamma(L^{\sigma^2})\dimm{3}\rto &
\Gamma(L)\dimm{3}\otimes\Gamma(L^\sigma\otimes L^{\sigma^2})\dimm{6}\rto & 
A_3\dimm{10}\rto & 0}}\end{equation}

By construction, $\Gamma\big((\Mm_PL)''\otimes L^\sigma\otimes
L^{\sigma^2}\big)\to \Gamma(L)\otimes \Gamma(L^{\sigma}\otimes
  L^{\sigma^2})$ factors through $W\otimes\Gamma(L^\sigma\otimes
  L^{\sigma^2})$.  We conclude:

\begin{prop}\Label{sadlemma}
If $P=Z(W)\in X$, then in $A_3$, we have $c_3\in WVV$.
\end{prop}

\subsubsection{Smooth point}

To go further, let us now assume that $P$ is a smooth point of
$X$. Then we have $\Mm_P L=L(-P)$, and $(\Mm_P
L)''=L^{-1}(P)$.  

Consider the  lattice of subspaces
$$\xymatrix{
\Gamma(L)\dimm{3}\otimes \Gamma(L^\sigma\otimes L^{\sigma^2})\dimm{6}\\
W\dimm{2}\otimes \Gamma(L^\sigma\otimes L^{\sigma^2})\dimm{6}\uto\\
\Gamma\big(L^{-1}(P)\otimes L^\sigma\otimes
L^{\sigma^2}\big)\dimm{4}\uto&W\dimm{2}\otimes
\Gamma\big(L^\sigma\otimes 
L^{\sigma^2}(-T)\big)\dimm{3}\ulto\\
\Gamma\big(L^{-1}(P)\otimes L^\sigma\big)\dimm{1}\otimes
\Gamma(L^{\sigma^2})\dimm{3}\uto^\alpha 
&\Gamma\big(L^{-1}(P)\otimes L^\sigma\otimes
L^{\sigma^2}(-T)\big)\dimm{1}\rlap{\,.}\uto\ulto}$$ 
Here $T$ is an effective  Cartier divisor  of degree 1 or~2. 
 
\begin{lem}
Let $Q\in X^\reg$ be the unique point such that $L^\sigma(P)=L(Q)$. Suppose
that $Q\not\in\supp T$, and that $L^{-1}(P)\otimes L^\sigma\otimes
L^{\sigma^2}(-T)$ is generated by global sections. Then we have 
\begin{multline*}
\Gamma\big(L^{-1}(P)\otimes L^\sigma\otimes
L^{\sigma^2}\big)\dimm{4}\\
=\,\Gamma\big(L^{-1}(P)\otimes 
L^\sigma\big)\dimm{1}\otimes\Gamma(L^{\sigma^2})\dimm{3}\,+
\,\Gamma\big(L^{-1}(P)\otimes  
L^\sigma\otimes L^{\sigma^2}(-T)\big)\dimm{1}\,.
\end{multline*}
\end{lem}
\begin{pf}
The line bundle $L^{-1}(P)\otimes L^\sigma$ has degree 1, so it has an essentially
unique non-zero section, which we shall call $g\in\Gamma\big(L^{-1}(P)\otimes
L^\sigma\big)$.  It vanishes  exactly at the  point $Q$. 
We conclude that $\Gamma\big(L^{-1}(P)\otimes
L^\sigma\big)=\Gamma\big(L^{-1}(P)\otimes 
L^\sigma(-Q)\big)$. Therefore the image of
$\Gamma\big(L^{-1}(P)\otimes
L^\sigma\big)\otimes\Gamma(L^{\sigma^2})$ in $\Gamma\big(L^{-1}(P)\otimes
L^\sigma\otimes L^{\sigma^2}\big)$ is 
contained in $\Gamma\big(L^{-1}(P)\otimes L^\sigma\otimes
L^{\sigma^2}(-Q)\big)$.  

As the bundle $L^{-1}(P)\otimes L^\sigma\otimes L^{\sigma^2}(-T)$
is generated by global sections, 
$\Gamma\big(L^{-1}(P)\otimes L^\sigma\otimes L^{\sigma^2}(-T)\big)$
is not contained in 
$\Gamma\big(L^{-1}(P)\otimes L^\sigma\otimes
L^{\sigma^2}(-Q)\big)$. This implies the 
lemma for dimension reasons. 
\end{pf}

\begin{cor}\Label{sameass}
Under the same assumptions, as 
subspaces of $W\dimm{2}\otimes {}\Gamma(L^\sigma\otimes
L^{\sigma^2})\dimm{6}$, we have 
$$\Gamma\big(L^{-1}(P)\otimes L^\sigma\otimes L^{\sigma^2}\big)\dimm{4}\,\subset\,
{}\Gamma\big(L^{-1}(P)\otimes L^\sigma\big)\dimm{1}\otimes {}\Gamma(L^{\sigma^2})\dimm{3}\,+\,{}W\dimm{2}\otimes {}\Gamma\big(L^\sigma\otimes
L^{\sigma^2}(-T)\big)\dimm{3}\,.$$
\end{cor}

\begin{cor}\Label{speccor}
If $P=Z(W)$ is a smooth point of $X$ such that $\sigma(P)=P$,  then in
$A_3$, we have $c_3\in WVW$.
\end{cor}
\begin{pf}
This follows from the above considerations, upon taking $T=P$, in
particular Sequence~(\ref{sec}), and the fact that $V\otimes W\to
\Gamma\big(L^\sigma\otimes L^{\sigma^2}(-P)\big)$ is surjective.
\end{pf}

\begin{cor}\Label{transnoin}
If $P=Z(W)$ is a smooth point of $X$, and $\sigma$ is a
translation, of order not equal to~2, then in $A_3$, we have $c_3\in
WWW$. 
\end{cor}
\begin{pf}
If $\sigma$ is the translation by the point $S\in X^\reg$, then
$Q=P-3S$, in the group law on $X^\reg$. Moreover, $P^\sigma=P-S$, and
$P^{\sigma^2}=P-2S$, so the lemma and its corollary apply to
$T=P+P^\sigma$, as $2S\not=0$. 

We also use that fact that $W\otimes W\to
\Gamma\big(L^\sigma(-P^\sigma)\otimes
L^{\sigma^2}(-P^{\sigma^2})\big)$ is surjective, which follows from
Lemma~\ref{themainlem}, below.
\end{pf}

\subsection{Non-special points}

This is the case where if $P\in X$, then $\sigma(P)\not=P$. We do not
exclude the case that $P$ is a singular point  of $X$.

\begin{lem}\Label{nonspw}
Let $W$ be a non-special 2-dimensional subspace of $V=B_1$. In $B_2$ we have
$WV+VW=VV=B_2$. 

Moreover, if $P\not\in X$, we have the stronger result that $B_2=WV$. 
\end{lem}
\begin{pf}
If the base locus of $W$ is empty, $W$ generates $L$, and we have a
short exact sequence of vector bundles 
on $X$: 
$$\xymatrix{
0\rto & M\rto & W\otimes \O_X\rto & L\rto & 0}\,.$$
Here $M$ is a line bundle which is isomorphic to $L^{-1}$, but not
canonically so. Tensoring with $L^\sigma$ and
taking global sections gives 
us the exact sequence
$$\xymatrix{
W\otimes \Gamma(L^\sigma)\rto & \Gamma(L\otimes L^\sigma)\rto & H^1(M\otimes
L^\sigma)}\,.$$
As $L\not\cong L^\sigma$, we have $M\otimes L^\sigma\not\cong\O_X$, and hence
$H^1(M\otimes L^\sigma)=0$, and so $W\otimes \Gamma(L^\sigma)\to
\Gamma(L\otimes L^\sigma)$ is surjective.  This proves that $WV=B_2$. (In
degenerate cases, this proof uses the fact that $L^{-1}\otimes L^\sigma$ is
{\em tame}, in the language of \cite{ATV}, which is proved in Lemma~7.18
of [ibid.].)

For a line bundle $M$ on $X$, and reduced point $Q\in X$, the set of
sections of $M$ vanishing at $Q$ is 
$\Gamma(X,M\otimes\Mm_Q)\subset\Gamma(X,M)$.  (We do not use the
notation $\Gamma\big(X,M(-Q)\big)$ to include the case where $Q$ is a
singularity, and hence not a Cartier divisor.) If $M$ has 
sections which do not vanish at $Q$, then
$\dim\Gamma(M\otimes\Mm_Q)=\dim\Gamma(M)-1$. 

Now suppose that $W$ has a base point $P\in X$. The base locus of
$W'=\sigma^\ast W\subset\Gamma(X,\sigma^\ast L)$ is
$P'=\sigma^{-1}P$. We have $W=\Gamma(L\otimes\Mm_P)\subset\Gamma(L)$ and
$W'=\Gamma(L^\sigma\otimes\Mm_{P'})\subset\Gamma(L^\sigma)$. Because $P'\not=P$,
the line bundle $L\otimes L^\sigma$ 
has sections which 
vanish at $P$ but not at $P'$, and sections which vanish at $P'$ but
not at $P$. Thus $\Gamma(L\otimes L^\sigma\otimes \Mm_P)$ and
$\Gamma(L\otimes L^\sigma\otimes \Mm_{P'})$ are two distinct codimension 1
subspaces of $\Gamma(L\otimes L^\sigma)$. Their sum is then necessarily the
whole space:
\begin{equation}\Label{cod1}
\Gamma(L\otimes \Mm_P\otimes L^\sigma)+\Gamma(L\otimes
L^\sigma\otimes\Mm_{P'})=\Gamma(L\otimes L^\sigma)\,.
\end{equation}
We now remark that the image of
$\Gamma(L\otimes\Mm_P)\otimes\Gamma(L^\sigma)$ in $\Gamma(L\otimes L^\sigma)$ 
is equal to $\Gamma(L\otimes\Mm_P\otimes L^\sigma)$. (This follows formally
from the fact that $\Gamma(L)\otimes \Gamma(L^\sigma)\to \Gamma(L\otimes
L^\sigma)$ is surjective, and that $\Gamma(L\otimes \Mm_P)\otimes
\Gamma(L^\sigma)\subset \Gamma(L)\otimes \Gamma(L^\sigma)$ and $\Gamma(L\otimes
L^\sigma\otimes \Mm_P)\subset\Gamma(L\otimes L^\sigma)$ are both codimension 1
subspaces.) Thus a reformulation of (\ref{cod1}) is
$WV'+VW'=VV'$. This is our claim.
\end{pf}

\subsubsection{An application to the study of \texorpdfstring{$c_3$}{c3}}

\begin{thm}\Label{anap}
Let $W\subset V$ be an arbitrary 2-dimensional subspace. 
In $A_3$, we have $c_3\in WWV+WVW+VWW$. 
\end{thm}
\begin{pf}
Let $P=Z(W)$.  If $P\not\in X$, we have $VV=WV$, and hence
$A_3=VVV=WVV=WWV$, and the claim is trivial. If $P\in X$, but
$\sigma(P)\not=P$, we use Proposition~\ref{sadlemma}, to deduce $c_3\in
WVV=W(WV+VW)$. This leaves the case where $P\in X$, and
$\sigma(P)=P$.  If $P$ is a non-singular point, our claim follows from
Corollary~\ref{speccor}.  If $P$ is a singular point, it follows
from~(\ref{crwww}). 
\end{pf}

We also need the following amplification:

\begin{prop}\Label{amplify}
Suppose $X$ is a Neron triangle, and $\sigma$ permutes the sides of
$X$ cyclically.  Assume that  $P=Z(W)$ is  a point on  $X$.  Then in $A_3$ we have
$c_3\in WWW$. 
\end{prop}
\begin{pf}
Fix a base point for $X$, making $X^\reg$ a group scheme, and $\sigma$
a translation.  Then the order of $\sigma$ is not equal to~2. So if
$P$ is a smooth point, the claim follows from
Corollary~\ref{transnoin}. If $P$ is a node, we can deduce the proof
from the case where $P$ is smooth by a degeneration
argument.\comment{may be fleshed out further}
\end{pf}

\section{Estimates}

Here we prove the stability estimates claimed in
Proposition~\ref{estimates}, by finding suitable values for the
functions $w_\ell(n)$ and $w_m(n)$, and then proving
Formula~(\ref{ggest}).

We may assume that $X$ is either smooth, or a Neron triangle, on which
$\sigma$ acts by cyclically permuting the sides. 

We always use notation $P=Z(W)$ and $Y=Z(U)$, and  treat the following
cases in turn: 
\begin{items}
\item $P$ lies off $X$,
\item $P$ lies on $X$, but is neither a node of $X$ nor a fixed point
  of $\sigma$, 
\item $P\in X^\reg$ is a fixed point of $\sigma$,
\item $P$ is a node of $X$,
\item $Y$ is a linear component of $X$.
\end{items}
In the cases (i) to (iv), we assume $Y\not\subset X$. Most of these
cases are divided into two subcases, a case where $\ell$ is large
compared to $m$, which includes the case $m=0$, and the converse case.

\subsection{The main dimension estimate}

\begin{lem}\Label{themainlem}
Let $L$ and $M$ be line bundles on $X$, both generated by global
sections, and both of degree at least~2. The multiplication map 
$$\Gamma(L)\otimes \Gamma(M)\longrightarrow \Gamma(L\otimes M)$$
is surjective, unless $\deg L=\deg M=2$ and $L\cong
M$. 
\end{lem}
\begin{pf}
Because of Proposition~\ref{717}, an equivalent statement is that 
$$H^1(L''\otimes M)=0\,.$$ The proof will proceed by induction on the
integer $\deg L+\deg M$.  The base case is when $\deg L=\deg M=2$. In
this case, $L''\cong L^{-1}$ is a tame line bundle. Then $L''\otimes M$ is also
tame. Hence $H^1(L''\otimes M)$ vanishes if and only if $L\not\cong M$. 

Now suppose that $\deg L>2$. Find a point $P\in X$, such that $L(-P)$
is still generated by global sections, and such that $L(-P)\not\cong
M$. We can apply the induction hypothesis to the bundles  $L(-P)$ and
$M$,  obtaining that 
$$H^1\big(L(-P)\otimes M''\big)=0\,.$$
But there is a surjection
$$\xymatrix{
H^1\big(L(-P)\otimes M''\big)\ar@{->>}[r]& 
H^1\big(L\otimes M''\big)}$$
proving that $H^1\big(L\otimes M''\big)=0$. 
\end{pf}

\begin{cor}\Label{themaindimension}
Let $L_1,\ldots,L_n$ be line bundles of degree at least 2 on $X$, all
generated by global sections.  The
map
$$\Gamma(L_1)\otimes\ldots\otimes \Gamma(L_n)\longrightarrow
\Gamma(L_1\otimes\ldots\otimes L_n)$$
is surjective, unless all $L_i$ are of degree 2 and pairwise
isomorphic to each other. 
\end{cor}
\begin{pf}
Assume that the bundles are not all of degree 2 or not all isomorphic
to one another. Then, after relabelling, we may assume that if $\deg
L_1=\deg L_2=2$, then $L_1\not\cong L_2$. Then we can apply the lemma 
successively to $L_1\otimes\ldots \otimes  L_{i}$ and $L_{i+1}$, for
$i=1,\ldots,n-1$. 
\end{pf}

\subsection{The test configuration of \texorpdfstring{$B$}{B}
  generated by a filtration}

We will write 
$$\{V^\gamma W^\alpha U^\beta\}\subset B_{\gamma+\alpha+\beta}$$
do denote the sum over all products $W_1\ldots
W_{\gamma+\alpha+\beta}$, where $\gamma$ of the $W_i$ are equal to
$V$, and $\alpha$ of the $W_i$ are equal to $W$, and $\beta$ of the
$W_i$ are equal to $U$. For example,
$$\{VWU\}=VWU+VUW+WUV+WVU+UVW+UWV\subset B_3\,.$$

\begin{lem}
Suppose that $\ell>0$ and $m>0$. Then we have
$$I^{(k)}_n=\sum_{\alpha,\beta} \{V^{n-\alpha-\beta}W^\alpha U^\beta\}\,,$$
where $(\alpha,\beta)$ ranges over all pairs of non-negative integers
satisfying
\begin{items}
\item $\alpha+\beta\leq n$,
\item \Label{bounds} $\alpha+\beta(\ell+1)\leq k\leq \alpha\ell+\beta(\ell+m)$.
\end{items}
\end{lem}
\begin{pf}
By construction, $\{V^{n-\alpha-\beta}V^\alpha U^\beta\}$ partakes in
$I_n^{(k)}$ if and only if there exist integers
$0<i_1,\ldots,i_\alpha\leq \ell$ and $\ell<j_1,\ldots,j_{\beta}\leq
\ell+m$ such that $k=i_1+\ldots+ i_\alpha+j_1+\ldots+j_\beta$. This
immediately proves the stated bounds on $k$. 

Conversely, if \ref{bounds} is satisfied, we can write $k=k_1+k_2$,
with $\alpha\leq k_1\leq \alpha\ell$ and $\beta(\ell+1)\leq k_2\leq
\beta(\ell+m)$. Then we can write $k_1=i_1+\ldots+i_\alpha$ with
$0<i_1,\ldots,i_\alpha\leq \ell$ and $k_2=j_1+\ldots+j_\beta$ with
$\ell<j_1,\ldots,j_{\beta}\leq 
\ell+m$.
\end{pf}

\begin{lem}\Label{alpha}
Suppose that $\ell>0$ and $m=0$. Then we have
$$I^{(k)}_n=\sum_{\alpha} \{V^{n-\alpha}W^\alpha\}\,,$$
where $\alpha$ ranges over all integers
such that $\alpha\leq n$ and 
$
\alpha\leq k\leq \alpha\ell$. \qed
\end{lem}

\begin{lem}\Label{beta}
Suppose that $\ell=0$ and $m>0$. Then we have
$$I^{(k)}_n=\sum_{\beta} \{V^{n-\beta}U^\beta\}\,,$$
where $\beta$ ranges over all integers
satisfying $\beta\leq n$ and $\beta\leq k\leq \beta m$. \qed
\end{lem}

\subsubsection{Useful reformulations: the \texorpdfstring{$W$}{W}-case}

\begin{lem}\Label{mleq2l}
For all $\alpha=1,\ldots,n$ we have
\begin{items}
\item if $(\alpha-1)\ell<k\leq \alpha\ell$ then
  $\{V^{n-\alpha}W^\alpha\}\subset I_n^{(k)}$, 
\item if $n\ell +(\alpha-1)m<k\leq  n\ell+\alpha m$ then
  $\{W^{n-\alpha}U^\alpha\}\subset 
  I_n^{(k)}$.
\end{items}
\end{lem}
\begin{pf}
The hardest claim is that $n\ell+(\alpha-1)m<k$ implies
$(n-i)+\alpha(\ell+1)\leq k$. This is proved by noting that under our
assumptions
$$0\leq (n-\alpha)(\ell-1)+(m-1)(\alpha-1)\,.$$
Reformulating gives
$$(n-\alpha)+\alpha(\ell+1)\leq n\ell+(\alpha-1)m+1\,,$$
which is what we need.
\end{pf}

The following variation will also be needed:

\begin{lem}\Label{positive}
For all $i=1,\ldots,n-1$, and every $k$, such that 
$$(i-1)\ell<k\leq i\ell\,,$$
we have $\{V^{n-i}W^i\}\subset I^{(k)}_n$.

For  all $i=1,\ldots, n-1$,  and all
$k$ which satisfy 
\begin{equation}\Label{assume}
(n-1)\ell+(i-1)m<  k\leq
(n-1)\ell+im\,\end{equation}
we have $\{V W^{n-1-i} U^{i}\}\subset I^{(k)}_n$. 

For $(n-1)(\ell+m)<k\leq (n-1)(\ell+m)+\ell$, we have $\{WU^{n-1}\}\subset
I_n^{(k)}$, and for $(n-1)(\ell+m)+\ell<k\leq n(\ell+m)$, we have
$U^n\subset I_n^{(k)}$. 
\end{lem}
\begin{pf}
Let us check the second claim. 
We have that
$$0\leq (\ell-1)(n-i-1)+(m-1)(i-1)\,,$$
because each of the four quantities involved is individually
non-negative, by our assumptions. Rewriting this inequality gives
\begin{equation}\Label{in1}
(n-1-i)+i(\ell+1)\leq(n-1)\ell+(i-1)m+1\,.
\end{equation}
We also have
\begin{equation}\Label{in2}
(n-1)\ell+im\leq (n-1-i)\ell+i(\ell+m)\,.
\end{equation}
Now suppose that $k$ satisfies (\ref{assume}). Then, in fact,
$$(n-1)\ell+(i-1)m+1\leq k\leq (n-1)\ell+im\,.$$
From (\ref{in1}) and (\ref{in2}) we conclude
$$(n-1-i)+i(\ell+1)\leq k\leq (n-1-i)\ell+i(\ell+m)\,.$$
With $\alpha=n-1-i$ and $\beta=i$ this is 
$$\alpha+\beta(\ell+1)\leq k\leq \alpha\ell+\beta(\ell+m)\,,$$
which is the condition we need.
\end{pf}

\subsubsection{The \texorpdfstring{$U$}{U}-case}

\begin{lem}\Label{mgeq2l}
Suppose that $m>0$. For all $i=1,\ldots,n$ we have
\begin{items}
\item if $(i-1)(\ell+m)<k\leq i(\ell+m)-m$ then $\{V^{n-i}WU^{i-1}\}\subset
  I_n^{(k)}$,
\item if $i(\ell+m)-m<k\leq i(\ell+m)$ then $\{V^{n-i}U^i\}\subset
  I_n^{(k)}$.
\end{items}
\end{lem}
\begin{pf}
Similar.
\end{pf}

\subsection{Exploiting non-commutativity}

Let  $Z(U)\not\subset X$.  The main dimension estimate
Corollary~\ref{themaindimension}, gives 
that 
\begin{align*}
\dim V^{n-\alpha-\beta} W^\alpha U^\beta
&\geq 3(n-\alpha-\beta)+2\alpha\\
&=3n-\alpha-3\beta\,,
\end{align*}
at least if $n-\alpha-\beta\geq1$. We can improve this estimate for
$\{V^{n-\alpha-\beta} W^\alpha U^\beta\}$.

\begin{prop}\Label{exploit}
We have 
\begin{items}
\item for $0\leq \beta\leq \frac{n}{2}$,
$$\dim \{ V^{n-\beta}U^\beta\}\geq 3n-2\beta\,,$$
\item for $\frac{n}{2}\leq \beta\leq n$,
$$\dim\{ V^{n-\beta}U^\beta \}\geq 4n-4\beta\,.$$
\item for $1\leq \beta\leq \frac{n+1}{2}$,
$$\dim \{ V^{n-\beta}WU^{\beta-1}\}\geq 3n+1-2\beta\,,$$
\item for $\frac{n+1}{2}\leq \beta\leq n$,
$$\dim\{ V^{n-\beta}WU^{\beta-1} \}\geq 4n+2-4\beta\,.$$
\end{items}
If $P=Z(W)$ is a smooth point of $X$, such that $\sigma (P)=P$, 
we also have 
\begin{items}
\item for $0\leq \beta\leq \frac{n}{2}$ 
$$\dim \{ W^{n-\beta} U^\beta\}\geq
 2n-\beta\,,$$ 
\item for $\frac{n}{2}\leq \beta\leq n$
$$\dim\{ W^{n-\beta} U^\beta \}\geq
  3n-3\beta\,.$$ 
\end{items}
\end{prop}
\begin{pf}
Choosing a non-zero element $s\in U$ defines an
injection $\O_X\to L$, which identifies $L$ with $\O(D)$, where $D$ is
the Cartier divisor on $X$, defined by the vanishing of $s$. 
Let $E$ be the greatest common divisor if $D$ and $D^\sigma$.  Because
$D\not=D^\sigma$, as $L\not\cong L^\sigma$, the degree of $E$ is at
most~2. Let 
$$H=\Gamma\big(X,L\otimes L^\sigma (-E)\big)\subset \Gamma(L\otimes
L^\sigma)\,.$$
Then an argument counting dimensions, using that
$UV=\Gamma\big(L\otimes L^\sigma(-D)\big)$ and
$VU=\Gamma(\big(L\otimes L^\sigma(-D^\sigma)\big)$,  proves $\dim
H=6-\deg E\geq 4$, and that that
inside $B_2$ we have  
$$H=UV+VU\,.$$

Now for $0\leq \beta\leq\frac{n}{2}$, we have 
$$\{V^{n-\beta}U^\beta\}\supset
\{V^{\beta}U^\beta\}V^{n-2\beta}\supset H^\beta V^{n-2\beta}\,,$$ 
and hence 
\begin{align*}
\dim \{V^{n-\beta}U^\beta\}
&\geq \dim  H^\beta V^{n-2\beta}\\
&\geq 3(n-2\beta)+ \beta\dim H\\
&\geq 3(n-2\beta)+4\beta\\
&=3n-2\beta\,,
\end{align*}
using Lemma~\ref{themaindimension}. 

Similarly, for $\frac{n}{2}\leq\beta\leq n$, we have
$$\{V^{n-\beta}U^\beta\}\supset
\{V^{n-\beta}U^{n-\beta}\}U^{2\beta-n}\supset H^{n-\beta} U^{2\beta-n}\,,$$ 
and hence 
\begin{align*}
\dim \{V^{n-\beta}U^\beta\}
&\geq \dim  H^{n-\beta} U^{2\beta-n}\\
&\geq 4(n-\beta)\\
&\geq 4n-4\beta\,.
\end{align*}

For the next two claims, involving a factor of $W$, we proceed
similarly. Let us explicate the case
$1\leq \beta\leq\frac{n+1}{2}$.  Here we have 
\begin{multline*}
\{V^{n-\beta}WU^{\beta-1}\}\supset
\{V^{n-\beta}U^{\beta-1}\}W\supset\\
\{V^{\beta-1}U^{\beta-1}\}V^{n+1-2\beta}W\supset H^{\beta-1}
V^{n+1-2\beta}W\,,
\end{multline*}
and hence 
\begin{align*}
\dim \{V^{n-\beta}WU^{\beta-1}\}
&\geq \dim  H^{\beta-1} V^{n+1-2\beta}W\\
&\geq 3(n+1-2\beta)+ (\beta-1)\dim H+2\\
&\geq 3n+5-6\beta+4(\beta-1)\\
&=3n+1-2\beta\,,
\end{align*}
using Lemma~\ref{themaindimension}, or in the case that $W$ is a node,
Lemma~\ref{borlem}, below.

The last two claims have the same proof as the first two, except for
we replace $L$ by $L(-P)$, throughout. This means we replace $D$ by $D-P$, and
the greatest common divisor $E$ of $D-P$ and $D^\sigma-P$ is of degree at
most~1. The line bundle $L(-1)\otimes L^\sigma(-P)(-E)$ has degree at
least~3, and therefore the proof goes through, mutatis mutandis.
\end{pf}

\subsection{Case: \tpdf{$P$}{P} off \tpdf{$X$}{X}}

This is the case where $P=Z(W)$ is not on $X$.  Then $Y=Z(U)$ cannot
be a component of $X$, so it intersects $X$ in a Cartier divisor. 

Let us first deal with the case where $m=0$.  By Lemma~\ref{mleq2l}, 
for every $\alpha=1,\ldots,n$, we have $\ell$ instances of $k$ where
$\{V^{n-\alpha}W^\alpha\}\subset I_n^{(k)}$.  For $\alpha<n$, we use
Lemma~\ref{nonspw} to conclude that $V^{n-\alpha}W^\alpha=V^n=B_n$, which has
dimension~$3n$.  For $\alpha=n$, the estimate  $\dim W^n\geq 2$, will be
sufficient.
Therefore, we may use
\begin{align}
w_\ell(n)&=3n(n-1)+2\Label{sameest}\\
&=3n^2-3n+2\,.\nonumber
\end{align}
Now by Theorem~\ref{anap}, we have that $c_3\in J_3^{(2\ell)}$, and so
we compute
\begin{align*}
(w_\ell+2\tilde a)(n)&=3n^2-3n+2+(n-1)(n-2)\\
&=4n^2-6n+4\,,
\end{align*}
and then 
\begin{align*}
G_\ell(n)
&=4\zeta_2(n)-6\zeta_1(n)+4\zeta_0(n)-n(n+1)(n+2)\\
&=4\zeta_2(n)-6\zeta_1(n)+4\zeta_0(n)-3\zeta_2(n)+3\zeta_1(n)-2\zeta_0(n)\\
&=\zeta_2(n)-3\zeta_1(n)+2\zeta_0(n)\,.
\end{align*}
This is $0$, for $n=2$, and positive for $n\geq3$, which is what we
needed to prove.

Now if $m>0$, we use Lemma~\ref{positive}. Then for every
$\alpha=1,\ldots,n-1$, we have $\ell$ 
instances of $k$ where $\{V^{n-\alpha}W^\alpha\}\subset I_n^{(k)}$. 
 We also have $\ell$
instances of $k$ where $WU^{n-1}\subset I_n^{(k)}$, which has
dimension~2. Thus, we get the same estimate (\ref{sameest}) we used
above, leading to the same value for $G_\ell$:
$$G_\ell(n)=
\zeta_2(n)-3\zeta_1(n)+2\zeta_0(n)\,.$$

We also get that for every $\beta=1,\ldots, n-1$, there are $m$
instances of $k$ where $\{VW^{n-1-\beta}U^\beta\}\subset
I_n^{(k)}$. Again by Lemma\ref{nonspw}, we have
$VW^{n-1-\beta}U^\beta=V^{n-\beta}U^\beta$, which has dimension
$3(n-\beta)$, by Lemma~\ref{themaindimension}. We add another $m$
instances of $U^n$, which has dimension~1.  This gives us
\begin{align*}
w_m(n)&=\sum_{\beta=1}^{n-1}(3n-3\beta)+1\\
&=\frac{3}{2}n^2-\frac{3}{2}n+1\,,
\end{align*}
and
\begin{align*}
G_m(n)&=\frac{3}{2}\zeta_2(n)-\frac{3}{2}\zeta_1(n)+\zeta_0(n)
-\frac{1}{2}n(n+1)(n+2)  \\
&=0\,.
\end{align*}

Together with the above, we get
\begin{align*}
\frac{1}{\ell}\big(\ell G_\ell(n)+ m G_m(n)\big)
&=G_\ell(n)+\frac{m}{\ell}G_m(n)\\
&=
\zeta_2(n)-3\zeta_1(n)+2\zeta_0(n)\,,
\end{align*}
which we have already remarked is  $0$, for $n=2$, and positive for
$n\geq3$. This finishes the case where $P$ is not on $X$.

\subsection{Case: \tpdf{$P$}{P} not a node, not a fixed point}

Now we assume that $P=Z(W)$ is on $X$, but is neither a singularity of
$X$, nor a fixed point of $\sigma$.  We also assume
that $Y=Z(U)$ is not contained in $X$, so it defines  a Cartier
divisor on $X$. 

We distinguish to subcases: $\ell\geq m\geq0$, and $m>\ell\geq0$. 

In the subcase $\ell\geq m\geq0$, the main fact we
use about $W$ is that $VV=VW+WV$ in 
$B_2=A_2$, see Lemma~\ref{nonspw}. We also use the main estimate
Lemma~\ref{themaindimension}.

For the subcase $m>\ell\geq 0$, we  exploit the
non-commutativity for $\{V^{n-\beta}U^\beta\}$ and 
$\{V^{n-\beta} WU^{\beta-1}\}$, using Lemma~\ref{exploit}.

Finally, it will be important that $c_3\in I_3^{(2\ell)}$.

\subsubsection{Subcase  \texorpdfstring{$\ell\geq m\geq0$}{l large}}

Let us first deal with the case where $m=0$. Then by
Lemma~\ref{mleq2l}, for each $\alpha=1,\ldots,n$, we have $\ell$
instances of $k$, such that $\{V^{n-\alpha}W^\alpha\} \subset
I_n^{(k)}$.  By Lemma~\ref{nonspw}, we have
\begin{items}
\item for $1\leq \alpha\leq \lfloor\frac{n}{2}\rfloor$, that
  $\{V^{n-\alpha}W^\alpha\}=V^n=B_n$, which has dimension $3n$,
\item for
 $\lfloor\frac{n}{2}\rfloor+1\leq\alpha\leq n-1$, that 
$\{V^{n-\alpha}W^\alpha\}\supset V^{2n-2\alpha} W^{2\alpha-n}$, and
  the latter has dimension at least $ 4n-2\alpha$, by
  Lemma~\ref{themaindimension}.
\end{items}
We also use that $\dim W^n\geq2$, which can be deduced, for example, from 
Lemma~\ref{themaindimension}, upon choice of a suitable $U$.
Altogether, we get
\begin{align}
w_\ell(n)&=\sum_{\alpha=1}^{\lfloor\frac{n}{2}\rfloor}3n\label{estabo}
+\sum_{\alpha=\lfloor\frac{n}{2}\rfloor+1}^{n-1}(4n-2\alpha)+2\\
&=\sum_{\alpha=1}^{n-1}(4n-2\alpha)-\nonumber
\sum_{\alpha=1}^{\lfloor\frac{n}{2}\rfloor}(n-2\alpha)+2\\  
&=4n(n-1)-n(n-1)-n\lfloor\sfrac{n}{2}\rfloor+\nonumber
\lfloor\sfrac{n}{2}\rfloor(\lfloor\sfrac{n}{2}\rfloor+1) + 2\\
&=\frac{11}{4}n^2-\frac{5}{2}n+2-\sfrac{1}{2}\{\sfrac{n}{2}\}\,.\nonumber
\end{align}
We have $c_3\in J_3^{(2\ell)}$, by Theorem~\ref{anap}, and so to
calculate $G_\ell$, we use
\begin{align*}
(w_\ell+2\tilde a)(n)&=
\frac{11}{4}n^2-\frac{5}{2}n+2-\sfrac{1}{2}\{\sfrac{n}{2}\}+(n-1)(n-2)\\
&=\frac{15}{4}n^2-\frac{11}{2}n+4-\sfrac{1}{2}\{\sfrac{n}{2}\}\,.
\end{align*}
This gives us 
\begin{align*}
G_\ell(n)&=\frac{15}{4}\zeta_2(n)-
\frac{11}{2}\zeta_1(n)+(4-\sfrac{1}{2}\{\sfrac{n}{2}\})
\zeta_0(n)-n(n+1)(n+2)\\
&=\frac{15}{4}\zeta_2(n)-
\frac{11}{2}\zeta_1(n)+(4-\sfrac{1}{2}\{\sfrac{n}{2}\})\zeta_0(n)
-3\zeta_2(n)+3\zeta_1(n)-2\zeta_0(n)\\
&=\frac{3}{4}\zeta_2(n)-
\frac{5}{2}\zeta_1(n)+(2-\sfrac{1}{2}\{\sfrac{n}{2}\})\zeta_0(n)\,.
\end{align*}
Thus $G_\ell(2)\geq0$, and $G_\ell(n)>0$, for $n\geq3$, which is what
we needed to prove.

In the case that $m>0$, we  use Lemma~\ref{positive}, instead of
Lemma~\ref{mleq2l}. This gives us,
for each   $\alpha=1,\ldots,n-1$,  altogether $\ell$
instances of $k$
where  $\{V^{n-\alpha}W^\alpha\} \subset I_n^{(k)}$, and another
$\ell$ instances of $k$ where $WU^{n-1}\subset I_n^{(k)}$.
Using
Lemma~\ref{nonspw}, as above, and the fact that $\dim WU^{n-1}=2$, we
get the same estimate (\ref{estabo}), which we used above, and hence
$$G_\ell(n)=\frac{3}{4}\zeta_2(n)-
\frac{5}{2}\zeta_1(n)+(2-\sfrac{1}{2}\{\sfrac{n}{2}\})\zeta_0(n)\,,$$
as above.

Lemma~\ref{positive} gives us, for each $\beta=1,\ldots n-1$,
also   $m$ instances of $k$ for which 
$\{VW^{n-1-\beta} U^\beta\}\subset
  I_n^{(k)}$.  For such $k$, we get 
$\dim I_n^{(k)}\geq 2n+2-2\beta$, for $1\leq \beta\leq n-2$ (because
of $VV=VW+WV$), and $3$,
for $\beta=n-1$.  We also have $m$ instances of $U^n$, for a
contribution of $1$. In total,
\begin{align*}
w_m(n)&\geq
\sum_{\beta=1}^{n-2}(2n+2-2\beta)+3+1\\
&=n^2+n-2\,.
\end{align*}
Hence
\begin{align*}
G_m(n)&\geq
\zeta_2(n)+\zeta_1(n) -2\zeta_0(n) -\sfrac{1}{2}n(n+1)(n+2)\\
&=\zeta_2(n)+\zeta_1(n) -2\zeta_0(n)
-\sfrac{3}{2}\zeta_2(n)+\sfrac{3}{2}\zeta_1(n)-\zeta_0(n)\\ 
&=-\frac{1}{2}\zeta_2(n)+\frac{5}{2}\zeta_1(n) -3\zeta_0(n)\,.
\end{align*}

If we now assume that $\ell\geq m$, then we get
\begin{align*}
\frac{1}{m}\big(\ell G_\ell(n)+m G_m(n)\big)
&\geq \frac{1}{m}\big(m G_\ell(n)+m G_m(n)\big)\\
&=G_\ell(n)+G_m(n)\\
&=\frac{1}{4}\,\zeta_2(n)-(1+\sfrac{1}{2}\{\sfrac{n}{2}\})\,\zeta_0(n)\,.
\end{align*}
This is equal to $0$ at $n=2$, but positive for $n\geq3$, which is
what we needed to prove.

\subsubsection{Subcase  \texorpdfstring{$m>\ell\geq0$}{m large}}

We use Lemma~\ref{mgeq2l}. This gives us, for every $\beta=1,\ldots,
n$, the existence of $\ell$ instances of $k$ where
$V^{n-\beta}WU^{\beta-1}\subset I_n^{(k)}$, and $m$ instances
of $k$ where $V^{n-\beta}U^\beta\subset I_n^{(k)}$.  We estimate these
dimensions with Lemma~\ref{exploit} to obtain

\begin{align}
w_\ell(n)&=
\sum_{\beta=1}^{\lfloor\frac{n+1}{2}\rfloor}(3n+1-2\beta) \label{startformula}
+\sum_{\beta=\lfloor\frac{n+1}{2}\rfloor+1}^n(4n+2-4\beta)\\
&=
\sum_{\beta=1}^{n}(4n+2-4\beta) \nonumber
-\sum_{\beta=1}^{\lfloor\frac{n+1}{2}\rfloor}(n+1-2\beta)\\
&=\nonumber
2n^2-\lfloor\sfrac{n+1}{2}\rfloor(\lceil\sfrac{n+1}{2}\rceil-1)\\
&= \frac{7}{4}n^2+\sfrac{1}{2}\{\sfrac{n}{2}\}\,.\nonumber
\end{align}
and (not forgetting $\dim U^n=1$)
\begin{align*}
w_m(n)&=
\sum_{\beta=1}^{\lfloor\frac{n}{2}\rfloor}(3n-2\beta)
+\sum_{\beta=\lfloor \frac{n}{2}\rfloor+1}^n(4n-4\beta)+1\\
&=\sum_{\beta=1}^n(4n-4\beta)
-\sum_{\beta=1}^{\lfloor\frac{n}{2}\rfloor}(n-2\beta)+1\\
&=2n^2-2n+1-\lfloor\sfrac{n}{2}\rfloor(\lceil\sfrac{n}{2}\rceil-1)\\
&=\frac{7}{4}n^2-\frac{3}{2}n+1-\sfrac{1}{2}\{\sfrac{n}{2}\}\,.
\end{align*}
We have 
\begin{align*}
(w_\ell+2\tilde a)(n)&=
\frac{7}{4}n^2+\sfrac{1}{2}\{\sfrac{n}{2}\}+(n-1)(n-2)\\
&=\frac{11}{4}n^2-3n+2+\sfrac{1}{2}\{\sfrac{n}{2}\}\,,
\end{align*}
which  gives 
\begin{align*}
G_\ell(n)
&=\frac{11}{4}\zeta_2(n)-3\zeta_1(n)+(2+\sfrac{1}{2}\{\sfrac{n}{2}\})\zeta_0(n)
-n(n+1)(n+2)\\
&=-\frac{1}{4}\zeta_2(n)+\sfrac{1}{2}\{\sfrac{n}{2}\}\zeta_0(n)\,.
\end{align*}

We also have
\begin{align*}
G_m(n)&=
\frac{7}{4}\zeta_2(n)-\frac{3}{2}\zeta_1(n)
+(1-\sfrac{1}{2}\{\sfrac{n}{2}\})\zeta_0(n)-\frac{1}{2}n(n+1)(n+2)\\
&=\frac{1}{4}\zeta_2(n)
-\sfrac{1}{2}\{\sfrac{n}{2}\}\zeta_0(n)\,.
\end{align*}

We observe that for any $n\geq2$, we have $G_m(n)>0$.  This proves the
case where $\ell=0$. If we now assume that $m>\ell$, then we have, for
$n\geq2$, 
$$\ell G_\ell(n)+m G_m(n)>\ell G_\ell(n)+\ell G_m(n)=0\,,$$
which proves what we need.

\subsection{Case: \tpdf{$P$}{P} smooth fixed point}

This is the case where $P=Z(W)\in X^\reg$, and $\sigma(P)=P$. 
In particular, $Z(U)\not\subset X$. 

The previous proof of the case $m>\ell\geq0$, which starts with
Formula~(\ref{startformula}), applies here, too, and so we only need to
deal with the case $\ell\geq m\geq0$.  This case is different than before, because
we cannot use $VV=VW+WV$ in the estimate for $w_\ell$, instead we
exploit non-commutativity using 
Lemma~\ref{exploit}, to improve the estimate on $w_m$.

So assume that $\ell\geq m\geq0$.

We use Lemma~\ref{mleq2l}. Thus we have, for every $\alpha=1,\ldots,n$,
a total of $\ell$ instances of $k$, where 
$\{V^{n-\alpha}W^\alpha\}\subset I_n^{(k)}$, and a total of $m$ instances
of $k$ where $\{W^{n-\alpha}U^\alpha\}\subset I_n^{(k)}$.

We have $\dim V^{n-\alpha}W^\alpha\geq3n-\alpha$, by
Lemma~\ref{themaindimension}, even for $\alpha=n$, because
the line bundles $L(-P)$ and $L^\sigma(-P)$ are not
isomorphic. Therefore,  we can take
\begin{align*}
w_\ell(n)&= \sum_{\alpha=1}^n (3n-\alpha)\\
&= \frac{5}{2}n^2-\frac{1}{2}n\,.
\end{align*}
Then 
\begin{align*}
(w_\ell+2\tilde a)(n) &=
\frac{5}{2}n^2-\frac{1}{2}n+(n-1)(n-2)\\
&=\frac{7}{2}n^2-\frac{7}{2}n+2\,,\end{align*}
and
\begin{align*}
G_\ell(n)&=
\frac{7}{2}\zeta_2(n)-\frac{7}{2}\zeta_1(n)+2\zeta_0(n)-n(n+1)(n+2)\\
&=\frac{1}{2}\zeta_2(n)-\frac{1}{2}\zeta_1(n)\,.
\end{align*}
We observe that this is positive, for all $n\geq2$, proving the $m=0$
case. 

Now assume that $m>0$. 
For $\beta\leq \frac{n}{2}$ we have $\dim\{W^{n-\beta}U^\beta\}\geq
2n-\beta$, and for $\beta\geq\frac{n}{2}$, we have
$\dim\{W^{n-\beta}U^\beta\}\geq 3n-3\beta$, by Lemma~\ref{exploit}. 
Therefore, (not forgetting that $\dim U^n=1$)
\begin{align*}
w_m(n)&=
\sum_{\beta=1}^{\lfloor \frac{n}{2}\rfloor}(2n-\beta)+
\sum_{\beta=\lfloor\frac{n}{2}\rfloor+1}^n(3n-3\beta)+1\\
&=
\sum_{\beta=1}^{n}(3n-3\beta)-\sum_{\beta=1}^{\lfloor
  \frac{n}{2}\rfloor}(n-2\beta)+1\\
& =
\frac{3}{2}n^2-\frac{3}{2}n+1-\lfloor\sfrac{n}{2}\rfloor (\lceil
\sfrac{n}{2}\rceil-1)\\
&=\frac{5}{4}n^2-n+1-\sfrac{1}{2}\{ \sfrac{n}{2}\}\,.
\end{align*}
Hence
\begin{align*}
G_m(n)
&=\frac{5}{4}\zeta_2(n)-\zeta_1(n)+ (1-\sfrac{1}{2}\{
\sfrac{n}{2}\})\zeta_0(n)
-\frac{1}{2}n(n+1)(n+2)\\
&=-\frac{1}{4}\zeta_2(n)+\frac{1}{2}\zeta_1(n)-\sfrac{1}{2}\{
\sfrac{n}{2}\}\,\zeta_0(n)\,.
\end{align*}
As we are in the case $\ell\geq m$, we get
\begin{align*}
\frac{1}{m}\big(\ell G_\ell(n)+m G_m(n)\big)
&\geq G_\ell(n)+G_m(n)\\
&=\frac{1}{4}\zeta_2(n)-\sfrac{1}{2}\{\sfrac{n}{2}\}\,\zeta_0(n)\,,
\end{align*}
which we already determined to be positive, for all $n\geq2$.

\subsection{Case: \tpdf{$P$}{P}  node}

Now let us deal with the case where $P=Z(W)$ is a node of $X$. 
We can still use the arguments from the above subcase $m>\ell\geq0$,
which starts with Equation~(\ref{startformula}), and consider this
case proved. 

So assume that $\ell\geq m\geq0$. 

We remark that $P$ is not fixed by $\sigma$, so we still have the convenient fact
$VV=VW+WV$. On the other hand, 
the estimate for the dimension of $V^{n-\alpha-\beta}
  W^{\alpha}U^\beta$ drops, but this is made up for by the fact
that $c_3\in I_3^{(3\ell)}$, by Proposition~\ref{amplify}. 

We use Lemma~\ref{mleq2l}. 
Using Lemma~\ref{borlem}, below, this gives the estimate\comment{have
  a forward reference}
\begin{align*}
\ol w_\ell(n)&= \sum_{\alpha=1}^{\lfloor\frac{n}{2}\rfloor} 3n
+\sum_{\alpha=\lfloor\frac{n}{2}\rfloor+1}^{n}(5n+1-4\alpha)\\
&=\sum_{\alpha=1}^{n}(5n+1-4\alpha)-
\sum_{\alpha=1}^{\lfloor\frac{n}{2}\rfloor} (2n+1-4\alpha)\\
&=3n^2-n-\lfloor\sfrac{n}{2}\rfloor(2\lceil\sfrac{n}{2}\rceil-1)\\
&=\frac{5}{2}n^2-\frac{1}{2}n\,.
\end{align*}
But note that this estimate is not optimal.  For
example, we counted the contribution of $W^n$ with $n+1$, but in fact,
it is larger than that, because we cannot have only one node
involved. 
It will be sufficient to improve $w_\ell(n)$ by adding 1 to it:
$$w_\ell(n)=\frac{5}{2}n^2-\frac{1}{2}n+1\,.$$
We have $c_3\in J_3^{(3\ell)}$,
and hence we compute
\begin{align*}
(w_\ell+3\tilde a)(n)
&=\frac{5}{2}n^2-\frac{1}{2}n+\frac{3}{2}(n-1)(n-2)+1\\
&=4n^2-5n+3+1\,.\\
\end{align*}
This gives us
\begin{align*}
G_\ell(n)
&=4\zeta_2(n)-5\zeta_1(n)+3\zeta_0(n)-n(n+1)(n+2)+\zeta_0(n)\\
&=\zeta_2(n)-2\zeta_1(n)+\zeta_0(n)+\zeta_0(n)\\
&=\frac{1}{6}\,n(2n^2+3n-3)-\{\sfrac{n}{3}\}+\zeta_0(n)\,.
\end{align*}
This is positive, for all $n\geq2$, proving the stability bound in the
case $m=0$.  

For the case $m>0$, we note that Lemma~\ref{mleq2l} also gives us 
\begin{align*}
w_m(n)
&=\sum_{\beta=1}^n(n+1-\beta)\\
&=\frac{1}{2}n^2+\frac{1}{2}n\,,
\end{align*}
and hence
\begin{align*}
G_m(n)
&=\frac{1}{2}\zeta_2(n)+\frac{1}{2}\zeta_1(n)-\frac{1}{2}n(n+1)(n+2)\\
&=-\zeta_2(n)+2\zeta_1(n)-\zeta_0(n)\,.
\end{align*}

Then we have
\begin{align*}
\frac{1}{m}\big(\ell G_\ell(n)+m G_m(n)\big)
&\geq G_\ell(n)+G_m(n)\\
&=\zeta_0(n)\\
&>0\,.
\end{align*}

\subsection{Case: \tpdf{$Y$}{Y} linear 
component of \tpdf{$X$}{X}}

We assume now that $Y=Z(U)\subset X$.  Then $Y$ is a linear component
of $X$, and $X$ is a Neron triangle. Moreover, the automorphism
$\sigma$ acts transitively on the three edges of $X$.  Any choice of
non-singular base point of $X$ turns $X^\reg$, the non-singular part
of $X$, into a commutative group scheme, such that $\sigma$ becomes a
translation. The point $P=Z(W)$ necessarily lies on $X$, but is not
fixed by $\sigma$ (of course, it may be a node).

To formulate dimension estimates, let $Y'=\sigma^{-1}Y$, and
$Y''=\sigma^{-2}Y$ be the other two edges of $X$. Let $Q=Y\cap Y'$,
$Q'=Y'\cap Y''$ and $Q''=Y''\cap Y$ be the three nodes of $X$.  Let
$M$ be a line bundle on $X$, whose degrees on the three edges,
$d,d',d''$, are all non-negative.  Denote by $\Gamma\big(M(-\beta
Y)\big)=\Gamma(M\otimes \iI_Y^\beta)\subset \Gamma(M)$ the subspace of
sections which vanish to order at least $\beta$ on $Y$, and
$\Gamma\big(M(-\alpha Q)\big)=\Gamma\big(M\otimes \Mm_Q^\alpha)\subset
\Gamma(M)$ the subspace of sections vanishing to order at least
$\alpha$ at $Q$. (The ideal sheaf of the closed subscheme $Y\subset X$
is locally principal, generated by $s$, for any non-zero element $s\in
U$.)

\begin{lem}\Label{borlem}
We have
\begin{items}
\item if $2\beta\leq d'+d''+1$, then $\dim \Gamma\big(M(-\beta
  Y)\big)=d'+d''+1-2\beta$,  
\item if $\beta+\gamma\leq d''+1$, then $\dim \Gamma\big(M(-\beta Y-\gamma
  Y')\big)=d''+1-\beta-\gamma$,
\item if $2\beta\leq d+d'+d''+1$, then $\dim \Gamma\big(M(-\beta
  Q)\big)=d+d'+d''+1-2\beta$, 
\item if $\beta+\gamma\leq d'+1$, and $\beta+\gamma\leq d''+1$, then
  $\dim \Gamma\big(M(-\beta Y-\gamma 
  Q')\big)=d'+d''+2-2\beta-2\gamma$,
\item if $\beta+\gamma\leq d'+1$, and $\beta+\gamma\leq d+d''+1$, then
  $\dim\Gamma\big(M(-\beta Q-\gamma
  Q')\big)=d+d'+d''+2-2\beta-2\gamma$, 
\item if $\alpha+\beta\leq d'+1$, and $\beta+\gamma\leq d''+1$, and
  $\alpha+\gamma\leq d+1$, then $\dim\Gamma\big(M(-\alpha Q-\beta Q'-\gamma
  Q'')\big)=d+d'+d''+3-2\alpha-2\beta-2\gamma$. 
\end{items}
In call cases, it is important, that $\alpha,\beta,\gamma\geq1$. 
\end{lem}
\begin{pf}
These formulas follow easily by breaking up $X$ into rational nodal
curves.
\end{pf}

\begin{cor}\Label{firhek}
For $1\leq\beta\leq \frac{n+1}{2}$, we have
$$\Gamma\big(L_n(-\beta Y)\big)+\Gamma\big(L_n(-\beta
Y')\big)=\Gamma\big(L_n(-\beta Q)\big)\,,$$
and 
$$\Gamma\big(L_n(-\beta Q)\big)+\Gamma\big(L_n(-\beta
Y'')\big)=\Gamma(L_n)\,,$$
and hence also
$$\Gamma\big(L_n(-\beta Y)\big)+\Gamma\big(L_n(-\beta Y')\big)
+\Gamma\big(L_n(-\beta Y'')\big)=\Gamma(L_n)\,.$$
\end{cor}

\begin{cor}\Label{stugjk}
If $1\leq \gamma\leq\alpha\leq\frac{n+1}{2}$, we have
$$\Gamma\big(L_n(-\alpha Y-\gamma Y')\big)+
\Gamma\big(L_n(-\alpha Y-\gamma Y'')\big)=
\Gamma\big(L_n(-\alpha Y-\gamma Q')\big)\,,$$
and 
$$\Gamma\big(L_n(-\alpha Y-\gamma Q')\big)+
\Gamma\big(L_n(-\alpha Y'-\gamma Q'')\big)=
\Gamma\big(L_n(-\alpha Q-\gamma Q'-\gamma Q'')\big)\,,$$
and also
$$\Gamma\big(L_n(-\alpha Q-\gamma Q'-\gamma Q'')\big)+
\Gamma\big(L_n(-\alpha Y''-\gamma Q)\big)=
\Gamma\big(L_n(-\gamma Q-\gamma Q'-\gamma Q'')\big)\,.$$
\end{cor}

\begin{cor}\Label{lastthree}
For $1\leq \beta\leq \lfloor\frac{n}{3}\rfloor$, we have that
$$\{V^{n-\beta}U^\beta\}=B_n\,.$$

For $\lfloor\frac{n}{3}\rfloor+1\leq \beta\leq
2\lfloor\frac{n}{3}\rfloor$, we have
$$\{V^{n-\beta}U^\beta\}\supset \Gamma\big(L_n(-iQ-iQ'-iQ'')\big)\,,$$
where $i=\beta-\lfloor\frac{n}{3}\rfloor$.

If $n\equiv 1\mod 3$, and $\beta=2\lfloor\frac{n}{3}\rfloor+1$, we
have
$$\{V^{n-\beta}U^\beta\}\supset
\Gamma\big(L_n(-(\lfloor\sfrac{n}{3}\rfloor+1)Y
-\lfloor\sfrac{n}{3}\rfloor Q')\big)\,,$$
if $n\equiv 2\mod 3$, and $\beta=2\lfloor\frac{n}{3}\rfloor+1$, we
have
\begin{equation}\Label{lath3}
\{V^{n-\beta}U^\beta\}\supset
\Gamma\big(L_n(-(\lfloor\sfrac{n}{3}\rfloor+1)Q
-\lfloor\sfrac{n}{3}\rfloor Q'-\lfloor\sfrac{n}{3}\rfloor Q'')\big)\,,
\end{equation}
and if $\beta=2\lfloor\frac{n}{3}\rfloor+2$, we
have
$$\{V^{n-\beta}U^\beta\}\supset
\Gamma\big(L_n(-(\lfloor\sfrac{n}{3}\rfloor+1)Y
-(\lfloor\sfrac{n}{3}\rfloor+1) Y')\big)\,.$$
\end{cor}
\begin{pf}
Note that 
$$\underbrace{(UVV)\ldots (UVV)}_{\beta}\,\underbrace{V\ldots
  V}_{n-3\beta}=  
\Gamma\big(L_n(-\beta Y)\big)\,,$$
and 
$$\underbrace{VUV\ldots VUV}_{\beta}\,\underbrace{V\ldots V}_{n-3\beta}
=\Gamma\big(L_n(-\beta Y')\big)\,,$$ 
and
$$\underbrace{VVU\ldots VVU}_{\beta}\,\underbrace{V\ldots V}_{n-3\beta}=
\Gamma\big(L_n(-\beta Y'')\big)\,.$$ 
The first claim now follows from Corollary~\ref{firhek}. 
 
The second claim follows from Corollary~\ref{stugjk} upon considering
$$\underbrace{UUV\ldots
  UUV}_{\beta-\lfloor\frac{n}{3}\rfloor}\,\underbrace{UVV\ldots
    UVV}_{2\lfloor\frac{n}{3} \rfloor-\beta} \,\underbrace{V\ldots
    V}_{n-3\lfloor\frac{n}{3}\rfloor } $$
and its 5 `cousins'. 

The last three claims follow by similar considerations. 
\end{pf}

\subsubsection{Subcase: \tpdf{$\ell=0$\,.}{l=0}}

To prove the stability estimates, let us start with the case that
$\ell=0$.  From Lemma~\ref{mgeq2l}, we obtain, for every
$\beta=1,\ldots,n$ exactly $m$ instances of $k$ where $\{V^{n-\beta}U^\beta\}\subset
I_n^{(k)}$. By Corollary~\ref{lastthree}, we have  (for
$n\geq1$):
\begin{align*}
w_m(n)&= \sum_{\beta=1}^{\lfloor \frac{n}{3}\rfloor} 3n
+\sum_{\beta=\lfloor \frac{n}{3}\rfloor+1}^{2\lfloor
  \frac{n}{3}\rfloor}
\big(3n+3-6(\beta-\lfloor\sfrac{n}{3}\rfloor)\big)\\
&=\sum_{\beta=1}^{\lfloor\frac{n}{3}\rfloor}
(6n+3-6\beta)\\
&=(6n+3)\lfloor \sfrac{n}{3}\rfloor-3\lfloor
\sfrac{n}{3}\rfloor(\lfloor \sfrac{n}{3}\rfloor+1) \\
&=\lfloor \sfrac{n}{3}\rfloor(6n-3\lfloor\sfrac{n}{3}\rfloor)\\
&=\frac{1}{3}(5n+3\{ \sfrac{n}{3}\})(n-3\{ \sfrac{n}{3}\})\\
&=\frac{1}{3}\cdot
\begin{cases}
{5}n^2&\text{if $n\equiv 0\mod 3$}\\
5n^2-4n-1&\text {if $n\equiv 1\mod 3$}\\
5n^2-8n-4&\text{if $n\equiv 2\mod 3$}\rlap{\,.}\end{cases}
\end{align*}
But we have not been optimal, yet, if $n\equiv 1\,(3)$ or $n\equiv
2\,(3)$. By Corollary~\ref{lastthree}, if $n\equiv 1\,(3)$, and $n\geq
4$, we can add $\frac{1}{3}(2n+4)$, 
but for $n=1$, we can only add~1 and not~2.  In order to have uniform
formulas, we will therefore add only $\frac{1}{3}(2n+1)$.

If $n\equiv 2\,(3)$, and $n\geq 5$, we can add $n+5$ and
$\frac{1}{3}(n+1)$, but if $n=2$, we can only add $5$ and $1$, and not
$7$ and $1$. Again, in the interest of uniform formulas, we add only
$n+3$ and $\frac{1}{3}(n+1)$, for a total of $\frac{2}{3}(2n+5)$.

In sum, we have
\begin{align*}
w_m(n)&=\frac{1}{3}\cdot
\begin{cases}
{5}n^2&\text{if $n\equiv 0\mod 3$}\\
5n^2-2n&\text {if $n\equiv 1\mod 3$}\\
5n^2-4n+6&\text{if $n\equiv 2\mod 3$}\rlap{\,.}\end{cases}
\end{align*}
Hence,
\begin{align*}
G_m(n)&=\frac{1}{3}\cdot
\begin{cases}
5\zeta_2(n)-\frac{3}{2}n(n+1)(n+2)&\text{if $n\equiv 0\mod 3$}\\
5\zeta_2(n)-2\zeta_1(n)-\frac{3}{2}n(n+1)(n+2)&\text
    {if $n\equiv 1\mod 3$}\\ 
5\zeta_2(n)-4\zeta_1(n)+6\zeta_0(n)-\frac{3}{2}n(n+1)(n+2)&\text{if
     $n\equiv 2\mod 3$}\rlap{\,.}\end{cases} \\
&=\frac{1}{6}
\begin{cases}
\zeta_2(n)+9\zeta_1(n)-6\zeta_0(n)&\text{if $n\equiv 0\mod 3$}\\
\zeta_2(n)+5\zeta_1(n)-6\zeta_0(n)&\text {if $n\equiv 1\mod 3$}\\
\zeta_2(n)+\zeta_1(n)+6\zeta_0(n)&\text{if $n\equiv 2\mod 3$}\rlap{\,.}\end{cases}
\end{align*}
A direct calculation, using the formulas of Lemma~\ref{sigmanottt},
shows that $m G_m(n)>0$, for all $n\geq2$.

\subsubsection{Subcase: \tpdf{$\ell>0$\,.}{l>0}}

Let us now deal with the case where $\ell>0$. Lemma~\ref{mgeq2l} gives
us, for every $\beta=0,\ldots,n-1$ precisely $m$ instances of $k$
where $\{V^{n-\beta-1}WU^{\beta}\}\subset I_n^{(k)}$. Using the two
facts that  $VV=WV+VW$ and $WUV+VUW=VUV$, we can replace the single
factor of $W$ by a factor of $V$, and still use the same arguments as
before, in the following cases:
\begin{items}
\item $0\leq \beta\leq \lfloor\frac{n}{3}\rfloor$, 
\item
  $\lfloor\frac{n}{3}\rfloor+1\leq\beta\leq2\lfloor\frac{n}{3}\rfloor-1$,
  or $n\equiv 2\mod 3$\,.
\end{items}

We get
\begin{align*}
w_\ell(n)&\geq \sum_{\beta=0}^{\lfloor \frac{n}{3}\rfloor} 3n
+\sum_{\beta=\lfloor \frac{n}{3}\rfloor+1}^{2\lfloor
  \frac{n}{3}\rfloor-1}
\big(3n+3-6(\beta-\lfloor\sfrac{n}{3}\rfloor)\big)
\\
&=\sum_{\beta=1}^{\lfloor \frac{n}{3}\rfloor}(6n+3-6\beta)+
3n-(3n+3-6\lfloor\sfrac{n}{3}\rfloor)\\
&=\lfloor\sfrac{n}{3}\rfloor(6n-3\lfloor\sfrac{n}{3}\rfloor)+
6\lfloor\sfrac{n}{3}\rfloor-3  \\
&=\frac{1}{3}\cdot
\begin{cases}
{5}n^2+6n-9&\text{if $n\equiv 0\mod 3$}\\
5n^2+2n-16&\text {if $n\equiv 1\mod 3$}\\
5n^2-2n-25&\text{if $n\equiv 2\mod 3$}\rlap{\,.}\end{cases}
\end{align*}
If $n\equiv2\mod3$, we can add a summand of
$3n+3-6\lfloor\frac{n}{3}\rfloor$, and we get instead
\begin{align*}
w_\ell(n)&=\frac{1}{3}\cdot
\begin{cases}
{5}n^2+6n-9&\text{if $n\equiv 0\mod 3$}\\
5n^2+2n-16&\text {if $n\equiv 1\mod 3$}\\
5n^2+n-4&\text{if $n\equiv 2\mod 3$}\rlap{\,.}\end{cases}
\end{align*}

These estimates are still not sufficient. 

For $n\equiv 1\,(3)$, we consider $\{V^{n-\beta-1}WU^\beta\}$ with
$\beta=2\lfloor\sfrac{n}{3}\rfloor$. It contains the three subspaces
$$\underbrace{UUV\ldots UUV}_{\lfloor\sfrac{n}{3}\rfloor}\,W\,,\qquad
\underbrace{UVU\ldots UVU}_{\lfloor\sfrac{n}{3}\rfloor}\,W\,,\qquad
\underbrace{VUU\ldots VUU}_{\lfloor\sfrac{n}{3}\rfloor}\,W\,.$$
Considering the first two, we see that the factor $W$ can be moved
into a position where it does not affect
the dimension calculation. This means that we get a dimension estimate
for the sum of these three spaces which is worse by 1, than the
dimension estimate $3n+3-6(\beta-\lfloor\sfrac{n}{3}\rfloor)$, which we used for
$\{V^{n-\beta}U^\beta\}$. We get an estimate of
$3n+2-6\lfloor\sfrac{n}{3}\rfloor=n+4$. 

For $n\equiv2\,(3)$, consider the contribution (\ref{lath3}), which
gave us a dimension estimate of $n+3$, for $n\geq2$.  In the presence
of a factor of $W$, this estimate drops by $2$ to $n+1$. 

This gives us the following improved formulas for $w_\ell$:
\begin{align*}
w_\ell(n)&=\frac{1}{3}\cdot
\begin{cases}
{5}n^2+6n-9&\text{if $n\equiv 0\mod 3$}\\
5n^2+5n-4&\text {if $n\equiv 1\mod 3$}\\
5n^2+4n-1&\text{if $n\equiv 2\mod 3$}\rlap{\,.}\end{cases}
\end{align*}
We have that $c_3\in J_3^{(3\ell)}$, by Proposition~\ref{amplify},
and so we add
$\frac{3}{2}(n-1)(n-2)$ to get
\begin{align*}
(w_\ell+3\tilde a)(n)&
=\frac{1}{6}\cdot
\begin{cases}
19n^2-15n&\text{if $n\equiv 0\mod 3$}\\
19n^2-17n+10&\text {if $n\equiv 1\mod 3$}\\
19n^2-19n+16&\text{if $n\equiv 2\mod 3$}\rlap{\,,}\end{cases}
\end{align*}
which gives, by subtracting $3\zeta_2-3\zeta_1+2\zeta_0$, 
\begin{align*}
G_\ell(n)&=\frac{1}{6}\cdot
\begin{cases}
\zeta_2(n)+3\zeta_1(n)-12\zeta_0(n)&\text{if $n\equiv 0\mod 3$}\\
\zeta_2(n)+\zeta_1(n)-2\zeta_0(n)&\text {if $n\equiv 1\mod 3$}\\
\zeta_2(n)-\zeta_1(n)+4\zeta_0(n)&\text{if $n\equiv 2\mod
  3$}\rlap{\,.}\end{cases} 
\end{align*}
A direct calculation shows that $\ell G_\ell(n)>0$, for all
$n\geq2$. This proves the required stability estimate for $m=0$.  

If both $\ell>0$, and $m>0$, Lemma~\ref{mgeq2l} shows that we may
simply add our values for $\ell G_\ell$ and $m G_m$ to prove the
stability estimate.  Since they are positive individually, this has
already been observed.


\section{Moduli stacks}

\subsection{Moduli stacks of stable regular algebras}

Recall that we are working over an algebraically closed field $\cc$,
of characteristic not equal to~2 or~3.

\begin{defn}
A {\bf flat family of stable regular algebras}, parametrized by the $\cc$-scheme
$T$,  is a graded vector bundle 
$\A=\bigoplus_{n=0}^\infty\A_n$ over $T$, endowed with the structure of
sheaf of graded $\O_X$-algebras, such that for every
$t\in T$, the fibre $\A_t$ is a stable regular algebra as in
Theorem~\ref{mnaefo}. 
\end{defn}

Let us denote the $\cc$-stack of families of stable regular algebras by
$\MM^{s,r}$.

\begin{thm}\Label{stackthm}
The  stack $\MM^{s,r}$ of flat families of stable regular algebras is a
smooth Deligne-Mumford stack of finite type. It has 4 components,
$\MM^{s,r}_A$, $\MM^{s,r}_B$, $\MM^{s,r}_E$, and $\MM^{s,r}_H$. Concretely,
\begin{description}
\item[(A)] $\MM^{s,r}_A=[U/G_{216}]$, where $U\subset\pp^2$ is the
\comment{can we say how $G_{216}$ acts on the lines and points removed?}
  complement of 4 concurrent lines, and  4 points  in $\pp^2$. The
  group  $G_{216}$ is the group of automorphisms of the oriented
  affine plane over $\ff_3$ (which has 216 elements).  It acts  $U$ via its quotient
  $SL_2(\ff_3)$, 
\item[(B)] $\MM^{s,r}_B=[V/\zz_4]$, where $V\subset \pp^1$ is the
  complement of 3 points in $\pp^1$, and   $\zz_4$ acts via its
  quotient $\zz_2$.  
\item[(E)] $\MM^{s,r}_E=B\zz_3$, 
\item[(H)] $\MM^{s,r}_H=B\zz_4$.
\end{description}
\end{thm}
\begin{pf}
Let $\A$ be a flat family of stable regular algebras parametrized by
the scheme $T$. As all members of $\A$ are elliptic, the triple $(X,
\sigma,L)$ of $\A$ is a flat family of 
elliptic triples. We have a short exact sequence of vector bundles on~$T$:
$$\xymatrix{
0\rto& R\rto& V\otimes V\rto & \pi_\ast(L\otimes L^\sigma)=\A_2\rto &
0}\,.$$
Here $\pi:X\to T$ is the structure morphism.  The algebra $\A$ is
recovered from $R\to V\otimes_{\O_T}V$ as the quotient of the tensor
algebra on $V$ divided by the two sided sheaf of ideals generated by
$R$. 

Now  $\sigma:X\to X$
induces a group automorphism\comment{$\Pic^0$ is line bundles of
  degree $0$ on each component}
$\Pic^0(\sigma):\Pic(X/T)^0\to\Pic(X/T)^0$. 
The order of $\Pic^0(\sigma)$ is necessarily finite, of order 1,2, 3,
or 4, as we are avoiding characteristic 2 and 3.
The order of $\Pic^0(\sigma)$ is also locally constant over $T$, and
so $T$ breaks up into 4 open and closed subschemes $T_i$, where
$T_i\subset T$ is the locus where the order of $\Pic^0(\sigma)$ is
$i$. This proves that $\MM^{r,s}$ also breaks up into 4 open and
closed substacks $\MM_i^{r,s}$. 

Let us deal with each of the 4 components in turn, and assume that the
order of $\Pic^0(\sigma)$ is constant. 

{\bf Case A. }
First, assume that the order of $\Pic^0(\sigma)$ is 1, so that
$\Pic^0(\sigma)$ is the identity of $\Pic^0(X/T)$. 
By our classification of stable algebras, the automorphism $\sigma$
acts transitively on the set of components in every fibre. Therefore,
by Theorem~II~3.2 in~\cite{DelRap}, for any
section $P$ of $X^\reg$, there exists
a unique structure of generalized elliptic curve on $X$, having $P$ as
origin, such that $\sigma$ acts as translation. In particular,
$X^\reg$ is a commutative group scheme, which acts on $X$, and there
is another section $S$ of $X^\reg$, such that $\sigma(Q)=Q+S$, for all
sections $Q$ of $X$. 

Using this, we prove that 
in the the fibered product 
$$\xymatrix{
Z\ar@{->}[rrr]\dto &&& T\dto^L\\
X\ar@{->}[rrr]^-{Q\mapsto \O(Q+\sigma Q+\sigma^{-1} Q)}&&&
\Pic^3(X)}$$
the scheme $Z$ is a form of the oriented affine plane over $\ff_3$. In
particular, $Z$ is  a finite \'etale cover $Z\to T$ of degree 9. 

Then, at least \'etale locally, we can assume that $P$ is a section of
$Z$, which we use to turn $X$ into a generalized elliptic curve.  then
$Z$ it equal to the scheme of 3-division points in $X$, in particular,
$Z$ is an oriented vector bundle over $\ff_3$. We can choose, \'etale
locally, an oriented basis $Q_1,Q_2$ for $Z$. 
We let $S$ be the section of $X^\reg$, such that $\sigma$ is
translation by $S$. From the fact that the algebra $\A$ is elliptic,
it follows that $S$ avoids $Z\subset X$. 

Thus, at least \'etale locally, we can associate to $\A$ a generalized
elliptic curve with full oriented level-3-structure, with an extra
point $S$ on it. The ambiguity is
in the choice of the oriented coordinate system $(P,Q_1,Q_2)$  for
the bundle of affine planes $Z$. 
The only restriction is that in singular fibres, $S$ avoids the
3-division sections, that   $S$
stays away from 
$X^\sing$, and that in the singular fibres,   $S$ avoids the 
connected component of $P$.

Thus, conversely, let Let $\MM(3)$ be the moduli scheme of generalized
elliptic curves with full oriented level-3-structure.  Let
$\EE(3)\to\MM(3)$ be the universal curve. It classifies quintuples
$(E,P,Q_1,Q_2,S)$, where $(E,P)$ is a generalized elliptic curve,
whose geometric fibres are smooth or triangles, $Q_1$ and $Q_2$ are
3-division points on $E$, forming an oriented bases for $E_3$, and
$S\in E$ is a further point on $E$. Let $\EE(3)^0\subset \EE(3)$ be
the open subscheme defined by the conditions that $S$ is not a node of
$E$, not in the component of identity in any triangle,  and not in
$E_3$. Using the Hesse family of elliptic curves, it is not hard to
identify $\EE(3)_0$ with the complement of 4 
concurrent lines and 4 points in $\pp^2$. 

The scheme $\EE(3)_0$ has the tautological triple $\big(E,\tau_S),
\O(3P)\big)$ over it. Associated to this triple is a flat family of
stable regular algebras $\A_{\EE(3)_0}$, parametrized by $\EE(3)_0$. 
The scheme of isomorphisms of $\A_{\EE(3)_0}$ is canonically
identified with the transformation groupoid of $G_{216}$ acting on
$\EE(3)_0$.\comment{make this explicit, esp. the action}

We have seen above, that the algebra $\A_{\EE(3)_0}$ is versal, i.e.,
every flat family of stable regular algebras is \'etale locally
induced from $\A_{\EE(3)_0}$.  This
finishes the proof that $\MM^{s,r}_1\cong[\EE(3)_0/G_{216}]$.

{\bf Case B.}  Let $\MM(2)_0$\comment{notation!}  be the moduli stack of
(non-singular) elliptic curves with a fixed 2-division point.  It
classifies triples $(E,P,Q)$, where $(E,P)$ is a smooth elliptic
curve, and $Q$ is a non-zero section of $E_2$.  Note that the Legendre
family of elliptic curves identifies $\MM(2)_0$ with $[V/\zz_4]$, where
$V\subset \pp^1$ is the complement of 3 points in
$\pp^1$.\comment{make this explicit}
To $(E,P,Q)$ we associate the elliptic triple $(X,\sigma,L)$, given by 
 $X=E$, $\sigma:X\to X$ is the involution
$R\mapsto Q-R$, whose fixed points are the 4 second roots of $Q$, and
$L=\O(2P+Q)$. The associated family of algebras is a flat family of
stable regular algebras of Type~B over $\MM(2)_0$. 

Conversely, let $\A\to T$ be a flat family of elliptic regular
algebras, with associated flat family of elliptic triples
$(X,\sigma,L)$. Assume that the order of $\Pic^0(\sigma)$ is 2.  Then
$\Pic^0(\sigma)$ is necessarily the multiplication by $-1$
automorphism. This implies that for any two local sections $Q,Q'$ of
$X$, the line bundles $\O_X(Q+\sigma(Q))$ and $\O_X(Q'+\sigma(Q'))$
are isomorphic. Thus there exists a unique, and therefore global,
section $P$, such that $\O(Q+\sigma(Q)+P)\cong L$, for any $Q$. This
makes $(X,P)$ into an elliptic curve, and the automorphism $\sigma$
into $R\mapsto \sigma(P)-R$. Regularity of the triple $(X,\sigma,L)$
implies that $Q=\sigma(P)$ is a 2-division point on $(X,P)$. We see
that $\MM^{s,r}_2\cong\MM(2)_0$. 

{\bf Cases E and H.} Left to the reader.
\end{pf}

\comment{it would be nice to explicitly write down the algebras in
  terms of coordinates $\langle a,b,c\rangle$ and $\lambda$. }

\subsection{Density}

\begin{prop}
Suppose that $\A$ is a flat family of graded $q$-truncated algebras, 
parametrized by the finite type $k$-scheme $T$.   Then the
locus of points $t\in T$, such that the fibre $\A_t$ is the
$q$-truncation of a stable  regular algebra
is open in $T$.
\end{prop}
\begin{pf}
By definition, 
$$\A=\bigoplus_{n=0}^q\A_q\,,$$
is a direct sum of vector bundles, where $\A_0=\O_T$. 
Let us assume that $t_0$ is a point of $T$, such that the fibre $\A|_{t_0}$
of $\A$ over $t_0$ is the truncation of a non-singular elliptic
regular algebra. We will prove that there exists an open neighbourhood $U$ of
$t_0$ in $T$, such that for every $t\in U$, the fibre $\A|_t$ is a
non-singular elliptic  regular algebra as well. 

The rank of the vector bundle $\A_n$ is a locally constant function on
$T$, and so by restricting to an open neighbourhood of $t_0$, we may assume
that it is constant.  It is then equal to $\frac{1}{2}(n+1)(n+2)$,
because it takes that value at $t_0$, by Formula~(1.15) in~\cite{AS}.

Let us write
$V=\A_1$. This is now a vector bundle of rank 3 on $T$.
Multiplication
in $\A$ defines a homomorphism of vector bundles
\begin{equation}\Label{mult2}
V\otimes_{\O_T}V\to \A_2\,.
\end{equation}
The  locus of points in $T$, over which (\ref{mult2}) is not
surjective is closed in $T$. Over the point $t_0$, the homomorphism
(\ref{mult2}) is surjective, because $\A|_{t_0}$ is generated in
degree 1.  So the locus of points in $T$ where (\ref{mult2}) is
surjective is an open neighbourhood of $t_0$, and be restricting to
this open neighbourhood, we may assume that (\ref{mult2}) is
surjective. Then the kernel of (\ref{mult2}) is a vector bundle
$R$ of rank~3 on $T$. 

By the same reasoning, we may assume that $\A_2\otimes_{\O_T} V\to \A_3$ 
and $V\otimes_{\O_T}\A_2\to\A_3$ are epimorphisms of vector bundles,
and that the respective kernels $K$ and $K'$ are vector bundles as
well.  Both vector bundles $K$ and $K'$ have rank~8.  

We have the following commutative diagram of coherent sheaves of
$\O_T$-modules:
$$\xymatrix{
W\rto^\beta \dto_\alpha & R\otimes V\rto\dto & K'\rto\dto &Q'\\
V\otimes R\rto\dto& V\otimes V\otimes V\rto\dto& V\otimes \A_2\dto\\
K\rto \dto & \A_2\otimes V\rto& \A_3\\
Q}$$
All tensor products are over $\O_T$. All rows and columns are exact,
when completed with zeros on all 
ends. By the snake lemma, $Q'=Q$.  Since $\A|_{t_0}$ is quadratic, the
homomorphisms $V\otimes R\to K$ and $R\otimes V\to K'$ are surjective
near $t_0$, and so we can assume they are surjective, and that
$Q'=Q=0$.  Then All sheaves in our diagram are locally free, and in
particular $W$ is a vector bundle of rank~1. 

We may also assume that the images of $\alpha$ and $\beta$ have full
rank, i.e., any non-zero section of $\omega$ of $W$ induces
isomorphisms $\beta(\omega)^\ast:R^\ast\to 
V$ and $\alpha(\omega)^\ast:V^\ast\to R$, because this is the case at
the point $t_0$.\comment{reference?}

We now consider the projective bundle $\pp(V)\to T$ of one-dimensional
quotients of $V$. From the homomorphism of vector bundles $R\to
V\otimes V$, we obtain homomorphisms
$$\Lambda^3R \longrightarrow \Lambda^3
V\otimes_{\O_T}\Sym^3
V\qquad\text{and}\qquad\Lambda^3R\longrightarrow\Sym^3
V\otimes_{\O_T}\Lambda^3V\,,$$ 
and hence
\begin{equation}\Label{lemfp}
\Lambda^3 V^\ast\otimes \Lambda^3R \longrightarrow \Sym^3
V\qquad\text{and}\qquad\Lambda^3R\otimes\Lambda^3 V^\ast\longrightarrow\Sym^3
V\,.\end{equation}
The arrows (\ref{lemfp}) are both strict monomorphisms of vector bundles
  over $T$, and hence $X_1$, $X_2$ are flat families of Cartier
  divisors of degree 3 in $\pp(V)$.
In fact, by our assumption that both $\alpha$ and $\beta$ have full
rank, we have $X_1=X_2$, i.e., we are in the semi-standard
case.\comment{I think this  
  is true, but I haven't proved it} We will call this scheme
$X=X_1=X_2$. 

Moreover, $R\to V\otimes V$ defines a family of subschemes $\Gamma$ 
in $\pp(V)\times_T\pp(V)$. By construction the projections
factor $\pi_1:\Gamma\to X$ and $\pi_2:\Gamma\to X$. At $t_0$, both
these morphisms are isomorphisms, and so by properness of $\Gamma$ and
$X$, we may assume that they are isomorphisms.

This gives us an
elliptic triple, from which we can construct a flat family of
elliptic regular algebras $\A'$, together with a morphism
of graded algebras $\A'\to \A$. At the point $t_0$ it is an
isomorphism, so we may assume that it is an isomorphism.

Let us finish by proving that there is an open neighbourhood of $t_0$
where the triple 
$(X,\sigma,L)$ is regular.  First note that the condition of being
regular or exceptional is open: it is the locus where $R^1\pi_\ast
(L\otimes (L^\sigma)^{-\otimes 2}\otimes L^{\sigma^2})$ has rank 1. 

If $\A_{t_0}$ is of type B, E, or H, we can pass to the open set where
$X$ is smooth, to get rid of the exceptional locus. 

So assume that $t_0$ is of type A. 
  If $X$ is smooth at $t_0$, we simply pass
to the open neighbourhood where $X$ is smooth.  If $X$ is not smooth at
$t_0$ it is a triangle, which is being rotated by $\sigma$. 
We now restrict to the open neighbourhood of $t_0$ where $X$ is nodal
and $\Pic^0(\sigma)$ is the identity. 

We also exclude all points where the algebra is not stable, as this
locus is closed. Then every fibre is either smooth, or a triangle on
which $\sigma$ acts by rotation, or exceptional, in which case
$\sigma$ acts by swapping components.  Thus in every fibre, $\sigma$ acts
transitively on connected components, and by by Lemme~II~1.7
of~\cite{DelRap}, we may endow $X$ with the structure of a generalized
elliptic curve.  Then by Proposition~II~1.15 of~\cite{DelRap}, the
exceptional locus, which is the locus where $X$ has two components, is
closed in $T$, and so we can remove it, and we are left with a family
of stable regular algebras, as required.
\end{pf}

\begin{cor}
For every $q$, the moduli stack of stable algebras $\MM_q^s$ has an
open substack isomorphic to 
$$\MM^{s,r}_A\amalg \MM^{s,r}_B\amalg\MM^{s,r}_E\amalg\MM^{s,r}_H\,.$$
Each stack 
$\MM^{s,r}_i$, for $i\in\{A,B,E,H\}$, is an open dense substack in the
irreducible component of $\MM_q^s$ which contains it. 
\end{cor}

We can unfortunately not prove that every stable algebra is
regular. There may be further stable algebras in the boundary of the
non-proper stacks $\MM^{s,r}_A$ and $\MM^{s,r}_B$, or there may be
entire components consisting of stable algebras which are not
regular. 

In particular, it may  happen that exceptional algebras are
stable, as their triples can satisfy our stability criterion, but we
have not examined this possibility closer. 

For every $q$, we have a projective coarse moduli space of
S-equivalence classes of semi-stable algebras, but we do not have a
complete description of these, except for low values of
$q$.\comment{see next section}


\end{document}